\documentclass[preprint,12pt]{elsarticle}

\usepackage{amssymb}
\usepackage{amsmath}
\usepackage{amsthm}
\newtheorem{theorem}{Proposition}
\usepackage{tikz}
\usepackage{float}
\usepackage{subcaption}
\journal{Nuclear Physics B}

\begin{document}

\begin{frontmatter}



\title{An Approximate Method for Modeling Plasmons on Graphene with Time- and Space-dependent Properties} 


\author[inst1]{Kaleigh Rudge}
\author[inst1]{Fadil Santosa} 
\address[inst1]{Department of Applied Mathematics and Statistics, Johns Hopkins University, Baltimore, MD 21218, USA}

\begin{abstract}
This paper presents approximate computational methods for modeling surface plasmon propagation on a graphene sheet with time- and space-dependent material properties. Starting from Maxwell's equations, the evolution of current density under a spatially and temporally modulated Drude weight is formulated as a partial integro-differential equation (PIDE). Assuming small-amplitude perturbations of the Drude weight, a regular perturbation expansion is developed. The solvability and Fourier invertibility of the resulting order-by-order integral equations are established in appropriate Sobolev spaces. A numerical solver combining a Volterra integral equation formulation with the fast Fourier transform (FFT) is introduced, along with pointwise error estimates and computational complexity analysis. Additionally, an analytical transform-based solution using contour integration and Laplace inversion is derived for traveling-wave Drude weight modulations. Numerical simulations demonstrate diverse wave phenomena, including traveling, standing, and exponentially growing plasmonic current densities.
\end{abstract}

    


\end{frontmatter}



\section{Introduction}
\label{intro} 
In this work, we develop an approximate method for calculating the current density on a graphene sheet associated with a propagating surface plasmon. The plasmon is manipulated by changing the Drude weight of the graphene sheet in time and in space. Our work is motivated by the need to develop efficient computational methods that could potentially serve as a tool for designing photonic devices that take advantage of plasmonic phenomena on 2D materials.

Interests in optical devices that leverage plasmonic effects on 2D material have been steadily growing in the past few years. The interest can be linked to the properties of surface plasmons on 2D systems such as graphene. They have been shown to be highly confined, low loss and tunable \cite{deabajo0,chen-et-al}. A recent work edited by Garcia de Abajo \cite{deabajo-et-al} describes the `road map' for graphene plasmonic devices. Among devices envisioned include sensors for bio applications and environmental monitoring, light emission, and microscopic imaging, made possible through integrated photonics.

One particularly exciting area is the development of integrative plasmonics which combines multi-physics \cite{zhou-et-al}. Some devices in this category are novel biosensors, energy conversion and photonic circuits. In \cite{esfandiari-et-al}, the authors put forward designs of graphene-based reconfigurable antennas, electro-magnetic wave absorbers, and sensors that exploit plasmonic waves. A multi-band perfect absorber is proposed in \cite{zeng-et-al}.

In the meantime, plasmonic effects have allowed for new designs of photodetectors and wave sources \cite{vitiello-viti}. The development of graphene plasmonic based biosensors and gas sensors continues \cite{wu-et-al, chi-et-al, tao-et-al, caridad-et-al,butt}.



The present work considers plasmon propagation on graphene whose material properties are modulated and continues the work in \cite{wilson-et-al, wilson-santosa-martin, santosa-shi}. Our interest is to see how the surface plasmons can be manipulated by perturbing the Drude weight associated with the 2D material. 

In \citep{wilson-santosa-martin}, an integro-differential equation is set up to model the modulation of the graphene when Drude weight is time-dependent. In \citep{santosa-shi}, the authors generalize the model to study the evolution of the current density on the graphene under a time- and space–dependent Drude weight. An integro-partial differential equation is derived for the current density. An explicit time-stepping finite difference method to solve an initial value problem for the resulting equation is proposed. 

We also consider a perturbational approach first introduce in \cite{tong-thesis}. The approach makes a major assumption which allows for simplification. Roughly speaking, it is assumed that the mixing of the different wave number components is limited. Numerical evidence provided some justification for the assumptions. In this work, we drop the assumption. Our work builds on both \citep{wilson-santosa-martin} and  \citep{santosa-shi}, solving the space and time-dependent integro-differential equation for the current density under a perturbation. 

The work in \cite{hu-et-al} is somewhat related to ours. The authors consider modulation of amplitude, wavelength and phase of surface plasmon polaritrons by changing the Fermi level. It should also be pointed out that the model problem under consideration can be thought of as a 1D wave propagation problem where the material property is time- and space-dependent. Analysis of such problems is somewhat sparse in open literature. \cite{wapenaar-aichele-vanmanen} provides an analysis that shows that time-only and space-only 1D wave equations are related, and therefore, tools for wave propagation in 1D inhomogeneous medium can be applied to study propagation problem where the material properties are time dependent.

This paper is organized as follows. We briefly explain our model in Section \ref{model_eqns}. In Section \ref{perturbation-expansion}, we derive the equations satisfied by the terms in a perturbation series. This is followed in Section \ref{solvability} by an analysis of the governing equations. Section \ref{discretization} describes the numerical method we propose to solve the perturbational equations, together with a rough assessment of the approximation errors. Section \ref{num-examples} contains numerical experiments which are similar to those considered in \cite{santosa-shi}. In Section \ref{special-form} we consider a special form of Drude weight perturbations which are in the form of a travelling wave. We show that the first order perturbation term can be solved using contour integrals and numerical quadrature. The main idea of the approach is described in the section. Some of the more technical calculations have been moved to Appendix A. Examples show that they recover the numerically computed solutions described in \cite{santosa-shi}. A short discussion section ends the paper.

\section{Model Equations}
\label{model_eqns}
As in \citep{santosa-shi}, we start with Maxwell's equation for the transverse magnetic (TM) mode in 2-D. The electric and magnetic fields have the form $E = (E_x,E_y,0)$ and $H = (0,0,H_z)$ respectively and are functions of time $t$ and space $(x,y)$. The electric permittivity $\epsilon$ and the magnetic permeability $\mu$ are constant. The graphene sheet lies on the entire $x$-axis. Away from the sheet, the electromagnetic field is governed by Maxwell's equation, and the sheet itself is modeled via the jump conditions along the $x$-axis. For $y\neq 0$ we have
\begin{subequations}
\begin{align}
    \mu\frac{\partial H_z}{\partial t} &= \frac{\partial E_x}{\partial y}-\frac{\partial E_y}{\partial x} ,\\
    \epsilon\frac{\partial E_x}{\partial t} &= \frac{\partial H_z}{\partial y} ,\\
    \epsilon \frac{\partial E_y}{\partial t} &= -\frac{\partial H_z}{\partial x} . 
\end{align}\label{Maxwell}
\end{subequations}
On the sheet, we have jump conditions
\begin{subequations}
\begin{align}
    [\![H_z]\!]_{y=0} &= j(x,t),\\
    [\![E_x]\!]_{y=0} &= 0, 
\end{align}\label{jump}
\end{subequations}
where $j$ is the current density. The current density satisfies Drude's law
\begin{equation}
    \frac{\partial j}{\partial t} = -\frac{1}{\tau}j + D(x,t)E_x(x,0,t), \label{current-density-eqn}
\end{equation}
where $D(x,t)$ is the Drude weight and $\tau$ is a damping factor. We emphasize that the Drude weight is both time- and space-dependent. 

The problem we wish to solve is an initial value problem for \eqref{Maxwell}, \eqref{jump} and \eqref{current-density-eqn}. The initial conditions for the electromagnetic fields and the current density are those associated with the plasmon solution when the Drude weight is constant $D_0$. Thus, the initial conditions are \cite{wilson-santosa-martin, santosa-shi}: 
\begin{subequations}\label{ic}
\begin{align}
    E_{x}(x,y,0) &= \frac{\gamma_0}{\epsilon s_0}e^{i\xi_0 x}e^{-\gamma_0 |y|} ,\\
    E_{y}(x,y,0) &= \frac{i\xi_0}{\epsilon s_0}\text{sgn}(y)e^{i\xi_0 x} e^{-\gamma_0 |y|} ,\\
    H_{z}(x,y,0) &= \text{sgn}(y)e^{i\xi_0 x}e^{-\gamma_0 |y|}\\
    j(x,0) &= 2e^{i\xi_0 x} . 
\end{align}
\end{subequations}
Here, $\xi_0$ is the wave number associated with the plasmon and is a problem parameter. The time-frequency $\omega_0=-i s_0$ of the constant Drude weight case is found by solving a quartic equation for $s_0$ (equation (3.8) in \cite{wilson-santosa-martin}); reproduced here
\[
s_0^4 - \frac{2}{\tau} s_0^3 + \left( \frac{1}{\tau^2}-\frac{\mu\epsilon D_0^2}{4}\right) s_0^2 - \frac{\xi_0^2 D_0^2}{4\epsilon^2} = 0,
\]
and the decay parameter $\gamma_0$ satisfies the dispersion relation (equation (3.4) in \cite{wilson-santosa-martin})
\[
\gamma_0^2 = \mu\epsilon s_0^2 + \xi_0^2.
\]
Our goal is to find the time dependent fields given $D(x,t)$ with $D(x,t\leq 0)=D_0$.

Motivated by the desire to come up with a computational method to solve the initial value problem that avoids truncation of 2-D space and implementation of absorbing boundary conditions, \cite{santosa-shi} derived a partial integro-differential equation (PIDE) for the current density on the graphene sheet. The derivation is presented in detail in \cite{santosa-shi}. It is shown that $j(x,t)$ satisfies
\begin{multline}
    \frac{\partial j}{\partial t} + \Big(\frac{1}{\tau} + \frac{\eta}{2}D(x,t)\Big)j(x,t) + \frac{\eta}{2}D(x,t)\mathcal{F}_{x}^{-1}\Big\{ k_1(\xi,t) *_t \Tilde{j}(\xi,t)\Big \} \\= \frac{\eta}{2}D(x,t)j_0(x,t) + \frac{\eta}{2}D(x,t)\Big(k_1(\xi_0,t)*_t j_0(x,t)\Big) + D(x,t)E_{x_0}(x,0,t) .\label{tong-eqn}
\end{multline}
In the third term on the left-hand side $\Tilde{j}(\xi,t)$ is the Fourier transform of $j(x,t)$. The inverse Fourier transform can be found in closed form as described in \cite{santosa-shi}. We leave it in the form above for this work. On the right-hand side of \eqref{tong-eqn}, $j_0$ and $E_{x_0}$ are the current density and electric field associated with the constant Drude weight $D_0$
\begin{subequations}
    \begin{align}
    E_{x_0}(x,0,t) &= \frac{\gamma_0}{\epsilon s_0}e^{i\xi_0 x}e^{-s_0 t}, \label{E0-eqn}\\
    j_0(x,t) &= 2e^{i\xi_0 x}e^{-s_0 t}\label{j0-eqn}.
\end{align}
\end{subequations}
The kernel $k_1(\xi,t)$ is defined as
\begin{equation*}
    k_1(\xi,t) = \int_0^t\frac{c\xi}{t'} J_1(c\xi t')dt' ,
\end{equation*}
where $J_1(z)$ is Bessel function of the first kind and $c=1/\sqrt{\mu\epsilon}$ is the speed of light in vacuum. 

Instead of solving the initial value problem for the electro-magnetic field and the current density through \eqref{Maxwell}, \eqref{jump} and \eqref{current-density-eqn}, and initial conditions \eqref{ic}, we solve \eqref{tong-eqn} for $j(x,t)$ with initial condition $j(x,0)=j_0(x,0)$. 

The main assumption in this work is that the Drude weight is of the form
\begin{equation}
    D(x,t) = D_0 +\alpha d(x,t) \label{D-eqn} 
\end{equation}
where $\alpha \ll 1$ and $\|d\|_{L^2([0,T]),H^1(\mathbb{R}))}$ is bounded. We will pursue a perturbative approach.

\section{Perturbation Expansion}
\label{perturbation-expansion}
We use the notation
\begin{equation*}
    \tilde{u}(\xi) = \int_{-\infty}^{\infty}u(x)e^{-ix\xi}dx,
\end{equation*}
We consider the case with no damping, i.e. let $\tau \rightarrow \infty$. Taking Fourier transform of
\eqref{tong-eqn}, we get
\begin{align}
    &\frac{\partial \tilde{j}}{\partial t} + \frac{\eta}{2}\tilde{D} *_{\xi} \tilde{j} + \frac{\eta}{2}\tilde{D} *_{\xi} \Big(k_1(\xi,t) *_t \tilde{j}(\xi,t)\Big) \label{model_ref_eqn} \\
    &= \frac{\eta}{2}\tilde{D} *_\xi \tilde{j_0} + \frac{\eta}{2}\tilde{D} *_{\xi}\Big(k_1(\xi_0,t)*_t\tilde{j_0}(\xi,t)\Big) + \tilde{D}(\xi,t)*_{\xi}\tilde{E}_{x_0}(\xi,0,t). 
   \nonumber 
\end{align}
Note that $\Tilde{D}(\xi,t) = D_0\delta(\xi) + \alpha \Tilde{d}(\xi,t)$ where $\delta(\xi)$ is the Dirac $\delta$-function.

We make the assumption that $\Tilde{j}(\xi,t)$ may be expressed as 
\begin{equation}
\Tilde{j}(\xi,t) = \Tilde{j_0}(\xi,t) + \alpha \Tilde{j_1}(\xi,t) + \alpha^2\Tilde{j_2}(\xi,t) + \cdots ,  
\label{ansatz}
\end{equation}
where $\alpha$ is the small parameter defined as in \eqref{D-eqn} and $\Tilde{j_0}$ is defined in \eqref{j0-eqn}. Then substituting our expression for  $\tilde{j}$ into \eqref{model_ref_eqn} and grouping like powers of $\alpha$, we obtain the following equations
\begin{subequations}

\begin{align}
\label{zeroth-ord-eqn}
&\mathcal{O}(1): \;\; \frac{\partial \Tilde{j_0}}{\partial t} + \frac{\eta D_0}{2}\Tilde{j_0}(\xi,t) + \frac{\eta D_0}{2}\Big(k_1(\xi,t) *_t \Tilde{j_0}(\xi,t)\Big) \nonumber \\ & \hspace{2cm} =\frac{\eta D_0}{2}\Tilde{j_0}(\xi,t) + \frac{\eta D_0}{2}\Big(k_1(\xi_0,t)*_t\Tilde{j_0}(\xi,t)\Big) +D_0\Tilde{E}_{x_0}(\xi,t), \\ \label{first-ord-eqn}
&\mathcal{O}(\alpha): \;\; \frac{\partial \Tilde{j_1}}{\partial t} + \frac{\eta D_0}{2}\Tilde{j_1}(\xi,t)  + \frac{\eta D_0}{2}k_1(\xi,t)*_t\Tilde{j_1}(\xi,t) \nonumber \\ &\hspace{2cm}=\Tilde{d}(\xi,t)*_{\xi}\Tilde{E}_{x_0}(\xi,0,t),\\
&\mathcal{O}(\alpha^n): \;\; \frac{\partial \Tilde{j_n}}{\partial t} + \frac{\eta D_0}{2}\Tilde{j_n}(\xi,t) + \frac{\eta D_0}{2}k_1(\xi,t)*_t\Tilde{j_n}(\xi,t) \nonumber \\&\hspace{2cm} = -\frac{\eta}{2}\Tilde{d}(\xi,t)*_\xi\Big(\Tilde{j}_{n-1}(\xi,t) + k_1(\xi,t)*_t \tilde{j}_{n-1}(\xi,t)\Big). \label{nth-ord-eqn}
\end{align}
\end{subequations}

The solution to \eqref{zeroth-ord-eqn} is simply the Fourier transform of \eqref{j0-eqn}
\[
\tilde{j}_0 = 2 \delta(\xi-\xi_0) e^{-s_0 t}.
\]
To solve the higher order equations, we proceed by reformulating each equation into an integral equation. First, we focus on the $\mathcal{O}(\alpha)$ equation, and then generalize the result to the higher-order equations. The equation for $\tilde{j}_1$ is rewritten as
\begin{equation*}
    \frac{\partial \tilde{j_1}}{\partial t} + \beta\tilde{j_1}(\xi,t)  = -\beta k_1(\xi,t)*_t\tilde{j}_1(\xi,t) +
    \tilde{d}(\xi,t)*_{\xi}\tilde{E}_{x_0}(\xi,0,t),
\end{equation*}
where $\beta = \frac{\eta D_0}{2}$. Viewing the right-hand side as a forcing term, we solve the first order ODE and arrive at
\begin{align*}
    \tilde{j}_1(\xi,t) = - \beta \int_0^t e^{-\beta(t-s)}(k_1(\xi,\cdot)*_t \tilde{j_1}(\xi,\cdot))(s)ds \\
    + \int_0^t e^{-\beta(t-s)}\tilde{d}(\xi,s)*_{\xi}\tilde{E}_{x_0}(\xi,0,s)ds.
\end{align*}
Defining the second term on the right-hand side as $L_0(\xi,t)$ we obtain the integral equation
\begin{equation*}
    \tilde{j}_1(\xi,t)  
    = - \int_0^t \left( \beta \int_{\tau}^t e^{-\beta(t-s)}k_1(\xi,s-\tau) ds \right) \tilde{j_1}(\xi,\tau) d\tau + L_0(\xi,t).
\end{equation*}
 This is an integral equation for $\tilde{j}_1$
\begin{equation}
    (I+K)\tilde{j}_1 = L_0 \label{j1-int-eqn}
\end{equation}
where $I$ is the identity mapping and $K$ is the integral operator 
\begin{equation}
    (K\tilde{j}_1)(\xi,t) = \int_0^t \left( \beta \int_{\tau}^t e^{-\beta(t-s)}k_1(\xi,s-\tau) ds \right) \tilde{j_1}(\xi,\tau) d\tau.\label{kern-def}
\end{equation}
This is essentially a Volterra Integral Equation \cite{brunner}. 

Next, we find the equations for $\tilde{j}_n(\xi,t)$ for $n>1$. From \eqref{nth-ord-eqn} we have the following
\begin{align*}
    \frac{\partial \tilde{j}_n}{\partial t} &+ \beta \tilde{j}_n(\xi,t) = - \beta k_1(\xi,t)*_t\tilde{j}_n(\xi,t)\\& -\frac{\eta}{2}\tilde{d}(\xi,t)*_{\xi}\Big(\tilde{j}_{n-1}(\xi,t) + k_1(\xi,t)*_t\tilde{j}_{n-1}(\xi,t)\Big).
\end{align*}
Define
\begin{equation*}
    L_{n-1}(\xi,t) :=-\int_0^t e^{-\beta(t-s)}\frac{\eta}{2}\tilde{d}(\xi,s)*_{\xi}\Big(\tilde{j}_{n-1}(\xi,s) + (k_1(\xi,\cdot)*_t\tilde{j}_{n-1}(\xi,\cdot))(s)\Big)ds.
\end{equation*}
This expression, which serves as a forcing term, is known if $\tilde{j}_{n-1}$ is known. Applying the same principle as for the $n=1$ case, we arrive at the integral equation for $\tilde{j}_n(\xi,t)$
\begin{equation}
    (I+K)\tilde{j}_n = L_{n-1} \label{int-eqn}
\end{equation}
where $I$ and $K$ are the same operators as defined above.
Hence, we have a way to recursively compute $\tilde{j}_n$ for any integer $n>1$. 

\section{Solution method}
\label{solvability}
The approach we propose here is to solve the integral equations \eqref{j1-int-eqn} and \eqref{int-eqn} for a given $\xi$ and invert the Fourier transform to obtain $j_n(x,t)$ for $n=1,2,\cdots$.

\subsection{Solvability of the integral equation}

Consider the integral equation \eqref{j1-int-eqn} on the interval $t\in [0,T]$. Existence and uniqueness of the solution for \eqref{j1-int-eqn} is established by showing that the kernel of $K$ is continuous in $[0,T]\times[0,T]$ and $L_0$ is continuous in $[0,T]$ \cite{brunner} (Theorem 1.2.3). The kernel in question is
\[
\beta e^{-\beta t} \int_\tau^t e^{\beta s} k_1(\xi,s-\tau) ds.
\]
It was shown in \citep{wilson-santosa-martin} that $k_1(\xi,t)$ is positive and continuous in $t$. 
It is now apparent that the kernel continuous in both $t$ and $\tau$. We note that if $\tilde{d}(\xi,t)$ is continuous in $t$, then $L_0(\xi,t)$ will be continuous, as $\tilde{E}_{x_0}(\xi,0,t)$ is exponential in $t$ and continuity is preserved under multiplication and integration. From these continuity results, we can conclude that we have a unique continuous solution $\tilde{j}_1(\xi,t)$ for $t \in [0,T]$.

To show the integral equation \eqref{int-eqn} is solvable, we need to show $L_{n-1}(\xi,t)$ is continuous in $t$. We have already established that $\tilde{j}_1(\xi,t)$ is continuous. Recall the exponential and $k_1$ are both continuous functions in $t$, so we can conclude that $L_1(\xi,t)$ is continuous. This in turn tells us that $\tilde{j}_2(\xi,t)$ is continuous. We apply this inductively to show that $L_{n-1}(\xi,t)$ is continuous. Therefore, we can conclude that $\tilde{j}_n(\xi,t)$ is unique and continuous for all integers $n \geq 1$.

The following proposition summarizes our result
\begin{theorem}
The integral equations \eqref{j1-int-eqn} and \eqref{int-eqn} admit unique continuous solutions $\tilde{j}_n(\xi,t)$ in $[0,T]$ for each $\xi$ for $\tilde{d}(\xi,t)$ that is continuous in $t$. 
\end{theorem}

\subsection{Invertibility of the Fourier transform}
We now show the inverse Fourier transform of $\tilde{j}_n(\xi,t)$ exists. We start by considering the invertibility of $\tilde{j}_1(\xi,t)$. The argument follows that in \cite{wilson-santosa-martin}. By multiplying \eqref{first-ord-eqn} by the complex conjugate of $\tilde{j}_1$ and adding the result to the equation obtained by multiplying the complex conjugate of \eqref{first-ord-eqn} by $\tilde{j}_1$, we get  
\begin{equation}
    \frac{d}{dt}|\tilde{j}_1|^2 + \eta D_0 |\tilde{j}_1|^2 + \eta D_0 \text{Re}[\Bar{\tilde{j}}_1\Big(k_1 *_t \tilde{j}_1\Big)] = 2\text{Re}[\Bar{\tilde{j}}_1\tilde{g}_0] ,
\end{equation}
where $\tilde{g}_0(\xi,t) = \tilde{d}(\xi,t)*_{\xi}\tilde{E}_{x_0}(\xi,t)$. Integrating in time over the interval $[0,t]$ leads to
\begin{equation}
    \frac{1}{2}|\tilde{j}_1(\xi,t)|^2 + \frac{\eta D_0}{2}\text{Re}[B_t(\tilde{j}_1,\tilde{j}_1)] = \text{Re}\langle\tilde{g}_0,\tilde{j}_1\rangle_{L^2(0,t)} , \label{en_eq_1}
\end{equation}
where 
\begin{equation*}
    B_t[u,v] = \langle u + k_1 *_t u, v \rangle_{L^2(0,t)} ,
\end{equation*}
and $\langle u, v \rangle_{L^2(0,t)} = \int_0^t u(s)\Bar{v}(s)ds$ \citep{wilson-santosa-martin}. 

To obtain an a priori estimate from \eqref{en_eq_1}, we need to show that the second term on the left-hand side is nonnegative. We rewrite the kernel $k_1$ as 
\begin{equation*}
    k_1(\xi,t) = c\xi\int_0^t \frac{J_1(c\xi s)}{s}ds = c\xi\int_0^{\infty}\frac{J_1(c\xi s)}{s}ds -c\xi\int_t^{\infty}\frac{J_1(c\xi s)}{s}ds.
\end{equation*}
The integral in first term is simply $\mathrm{sgn}(c\xi)$. Therefore, we can write
\[
k_1(\xi,t) = c|\xi| - c\xi\int_t^{\infty}\frac{J_1(c\xi s)}{s}ds .
\]
Now notice the second term in is even with respect to $\xi$. We can further rewrite
\[
k_1(\xi,t) = c|\xi| - c |\xi|\int_t^{\infty}\frac{J_1(c|\xi| s)}{s}ds .
\]
We apply Lemma 5.7 in \citep{wilson-santosa-martin} with $a=c|\xi|$ to conclude that the second term on the left-hand side of \eqref{en_eq_1} is nonnegative. So, from \eqref{en_eq_1} we have 
\[
\frac{1}{2}|\tilde{j}_1(\xi,t)|^2 \leq \text{Re}\langle \tilde{g}_0,\tilde{j}_1 \rangle_{L^2(0,t)} .
\]
To obtain an a priori estimate, we bound the right-hand side using Cauchy-Schwarz
\begin{align*}
\frac{1}{2}|\tilde{j}_1(\xi,t)|^2 &\leq \|\tilde{g}_0\|_{L^2(0,t)} \|\tilde{j}_1\|_{L^2(0,t)} \\
& \leq \sigma \|\tilde{g}_0\|^2_{L^2(0,t)} + \frac{1}{\sigma} \int_0^t |\tilde{j}_1(\xi,s)|^2 ds,
\end{align*}
for some $\sigma > 0$. Now, apply Gronwall's inequality to obtain the following bound
\begin{equation}\label{boundj1}
    |\tilde{j}_1(\xi,t)|^2 \leq 2\sigma\|\tilde{g}_0 \|^2_{L^2(0,t)} e^{2t/\sigma}.
\end{equation}
Next, integrate both sides of the inequality in $\xi$
\begin{equation}\label{l2-boundj1}
    \int_{-\infty}^\infty  |\tilde{j}_1(\xi,t)|^2 d\xi \leq 2\sigma e^{2t/\sigma} \int_{-\infty}^\infty \|\tilde{g}_0(\xi,\cdot)\|^2_{L^2(0,t)} d\xi .
\end{equation}
We conclude then that
\begin{equation}
    \|j_1(\cdot,t)\|^2_{L^2(\mathbb{R})} \leq 2\sigma e^{2t/\sigma} \|g_0\|^2_{L^2((0,t),L^2(\mathbb{R}))}. \label{apriori}
\end{equation}
Recall $g_0(x,t) = d(x,t)E_{x_0}(x,0,t)$. Therefore,
\begin{align*}
    \|g_0\|^2_{L^2((0,t),L^2(\mathbb{R}))} &= \int_0^t \int_{-\infty}^\infty |d(x,s)|^2 |E_{x_0}(x,s)|^2 dx ds \\
    &\leq \left|\frac{\gamma_0}{\epsilon\omega_0}\right|^2 \| d \|^2_{L^2((0,t),L^2(\mathbb{R}))},
\end{align*}
using \eqref{E0-eqn}. If $d$ is in $L^2((0,t),L^2(\mathbb{R}))$, by \eqref{apriori}, we arrive at the conclusion that the inverse Fourier transform of $\tilde{j}_1(\xi,t)$ exists for each finite $t$. Moreover, we can integrate \eqref{apriori} to assert that $\|j_1\|_{L^2((0,t),L^2(\mathbb{R}))} < \infty$.

Regularity of $j_1$ in space can be established if we make spatial smoothness assumption on $d$. From \eqref{boundj1}, we have
\[
\int_{-\infty}^\infty \xi^2 |\tilde{j}_1(\xi,t)|^2 d\xi \leq 2\sigma e^{2t/\sigma}\int_{-\infty}^\infty \xi^2 \| \tilde{g}_0(\xi,\cdot)\|^2_{L^2(0,t)} d\xi.
\]
Using $\tilde{g}_0(\xi,t)=\tilde{d(\xi,t)}*_\xi \tilde{E}_{x_0}(\xi,0,t)$ and $\tilde{E}_{x_0}(\xi,0,t)=\frac{\gamma_0}{\epsilon s_0}\delta(\xi-\xi_0)e^{-s_0 t}$, we get 
\[
\int_{-\infty}^\infty \xi^2 |\tilde{j}_1(\xi,t)|^2 d\xi \leq 
\int_{-\infty}^\infty \xi^2 \left| \tilde{d}(\xi-\xi_0)\frac{\gamma_0}{\epsilon s_0} e^{-s_0 t} \right|^2 d\xi.
\]
From this, we conclude that
\begin{equation}
  \int_{-\infty}^\infty \xi^2 |\tilde{j}_1(\xi,t)|^2 d\xi \leq
  \frac{\gamma_0^2}{\epsilon^2|s_0|^3} \int_{-\infty}^\infty  \xi^2|\tilde{d}(\xi-\xi_0,t)|^2 d\xi .
  \label{boundj1-H1}
\end{equation}

Next, we show that $\tilde{j}_n$ for $n>1$ are invertible using induction. First, consider \eqref{nth-ord-eqn} with $n=2$ as the base case. Observe the left hand side of \eqref{nth-ord-eqn} is of the same structure as the \eqref{first-ord-eqn}, so the same conjugation and multiplication manipulation can be made to arrive at
\begin{equation}
    \frac{1}{2}|\tilde{j}_2(\xi,t)|^2 + \frac{\eta D_0}{2} \text{Re}[B_t(\tilde{j}_2,\tilde{j}_2)] = \frac{-\eta}{2}\text{Re}\langle\tilde{d}*_\xi (\tilde{j}_1 + k_1 *_t \tilde{j}_1),\tilde{j}_2 \rangle_{L^2(0,t)}.
\end{equation}
Let $\tilde{g}_1 = \tilde{d} *_\xi (\tilde{j}_1 + k_1*_t\tilde{j}_1)$. Then using the same argument as for $\tilde{j}_1$, we are led to
\begin{equation} \label{boundj2}
    |\tilde{j}_2(\xi,t)|^2 \leq 2\sigma\|\tilde{g}_1 \|^2_{L^2(0,t)} e^{2t/\sigma}.
\end{equation}
What remains is to show that $\|g_1\|_{L^2((0,t),L^2(\mathbb{R}))}$ is bounded.
 
Taking the inverse Fourier transform of $\tilde{g}_1$, we have
\begin{equation}\label{g1-def}
g_1(x,t) = d(x,t)\Big(j_1(x,t) + \mathcal{F}^{-1}_{x}[k_1(\xi,\cdot)*_t \tilde{j}_1(\xi,\cdot)](x,t)\Big).
\end{equation}
Taking the norm of both sides, we have
\begin{align*} 
\|g_1\|_{L^2((0,t),L^2(\mathbb{R}))} 
&\leq \|d \; j_1\|_{L^2((0,t),L^2(\mathbb{R}))} \\ &\hspace{60pt}+ \|d \; \mathcal{F}^{-1}_{x}[k_1(\xi,\cdot)*_t \tilde{j}_1(\xi,\cdot)](\cdot,\cdot)\|_{L^2((0,t),L^2(\mathbb{R}))}.
\end{align*}
The first term on the right-hand side is bounded since both $d(x,t)$ and $j_1(x,t)$ have bounded $L^2$ norms. For the second term, we need to show that
\begin{align}
    \| \mathcal{F}^{-1}_{x}[k_1(\xi,\cdot)*_t & \tilde{j}_1(\xi,\cdot)](\cdot,\cdot)\|_{L^2((0,t),L^2(\mathbb{R}))}^2 \nonumber \\ & = \int_{-\infty}^\infty \int_0^t \left| 
    k_1(\xi,\cdot)*_t \tilde{j}_1(\xi,\cdot)(t') \right|^2 dt' d\xi, \label{toshow}
\end{align}
is bounded. 

We start with the expression for $k_1$ and consider
\begin{align*}
    \tilde{j}_1(\xi,\cdot)*_t& k_1(\xi,\cdot)(t)  = \int_0^t \tilde{j}_1(\xi,t-s) \left( c\xi\int_0^s\frac{J_1(c\xi\tau)}{\tau}d\tau \right) ds\\
    & =\int_0^t \tilde{j}_1(\xi,t-s) \left( c\xi \sum_{m=0}^{\infty}\frac{(-1)^m}{m!(m+1)!}\Big(\frac{c\xi}{2}\Big)^{2m+1}\frac{s^{2m+1}}{2m+1} \right)ds,
\end{align*}
after using the series expansion of $J_1(x)$ \cite{bateman} and integrating in $\tau$. Therefore, we have
\begin{align*}
& \left|   \tilde{j}_1(\xi,\cdot)*_t k_1(\xi,\cdot)(t)  \right| \\ 
    &\leq \int_0^t \Big | \tilde{j}_1(\xi,t-s)\Big| \Big | c\xi\sum_{m=0}^{\infty}\frac{(-1)^m}{m!(m+1)!}\Big(\frac{c\xi}{2}\Big)^{2m+1}\frac{s^{2m+1}}{2m+1}\Big | ds\\
        &\leq \int_0^t \Big | \tilde{j}_1(\xi,t-s)\Big| \Big | \frac{c^2\xi^2 s}{2}\sum_{m=0}^{\infty}\frac{1}{m!}\Big(\frac{-c^2\xi^2 s^2}{2}\Big)^m\frac{1}{(m+1)!(2m+1)}\Big | ds\\
    &\leq \int_0^t \Big| \tilde{j}_1(\xi,t-s)\Big| \Big| \frac{c^2\xi^2 s}{2}\sum_{m=0}^{\infty}\frac{1}{m!}\Big(\frac{-c^2\xi^2 s^2}{2}\Big)^m\Big|ds\\
    &= \int_0^t \Big|\tilde{j}_1(\xi,t-s)\Big|\Big| \frac{c^2\xi^2 s}{2}e^{-(c\xi s)^2/2}\Big |ds.
\end{align*}
A straightforward calculation shows that
\[
\frac{c^2\xi^2 s}{2}e^{-(c\xi s)^2/2} \leq 
\frac{c\xi}{2}e^{-1/2}, \;\; 0<s<t.
\]
Therefore, using Cauchy-Schwarz, we get
\begin{equation*}
    \left|   \tilde{j}_1(\xi,\cdot)*_t k_1(\xi,\cdot)(t)  \right| \leq
    \left\{ \frac{c^2}{4e} \int_0^t \xi^2 \left| \tilde{j}_1(\xi,t)\right|^2 dt \right\}^{1/2}.
\end{equation*}
Applying this in \eqref{toshow}, we have
\begin{align*} 
\| \mathcal{F}^{-1}_{x}[k_1(\xi,\cdot) &*_t \tilde{j}_1(\xi,\cdot)](\cdot,\cdot)\|_{L^2((0,t),L^2(\mathbb{R}))}^2\\
& \leq  \frac{c^2}{4e} \int_0^t \int_0^{t'} \left[ \int_{-\infty}^\infty \xi^2 |\tilde{j}_1(\xi,s)|^2 d\xi \right] ds dt'  \\
& \leq \frac{c^2}{4e} \frac{\gamma_0^2}{\epsilon^2|s_0|^3} \int_0^t \int_0^{t'} \left[ \int_{- \infty}^\infty\xi^2|\tilde{d}(\xi-\xi_0,s)|^2 d\xi \right] ds dt',
\end{align*}
by \eqref{boundj1-H1}. From this, we can conclude that if $d(x,t)$ is in $L^2((0,t),H^1(\mathbb{R}))$, then $g_1(x,t)$, defined in \eqref{g1-def}, is in $L^2((0,t),L^2(\mathbb{R}))$. Therefore, $\tilde{j}_2(\xi,t)$ is invertible.

This argument can be inductively applied to higher order terms to show the norm of $\tilde{j}_n(\xi,t)$ is invertible. We summarize our findings in
\begin{theorem}
For Drude weight perturbation $d(x,t)$ in $L^2((0,t),H^1(\mathbb{R}))$, the perturbational corrections $\tilde{j}_n(\xi,t)$ for $n=1,2,\dots$ have Fourier inverses for finite $t$ and are in $L^2((0,t),L^2(\mathbb{R}))$.
\end{theorem}

\section{Discretization}
\label{discretization}

Following the solution method described in the previous section, we will find approximate solutions to the integral equations \eqref{j1-int-eqn} and \eqref{int-eqn}. We discretize time $t$ and turn the integral operator defined in \eqref{kern-def} into a matrix-vector product. The time domain is discretized by fixing the time increment $\Delta t$ and the number steps. We will obtain approximate values of $\tilde{j}_n(\xi,t_m)$ where $t_m=m\Delta t$ for $m=1,\cdots,M$. Solving the integral equations amount to inverting a matrix.

Next, the spatial range is fixed -- we choose $[-4\pi,4\pi]$. Then, the number of points, denoted as $N$, is chosen. This determines the spatial increment $\Delta x=8\pi/N$, and the sample points are $x_l=-4\pi+(l-1)\Delta x$, $l=1,\cdots,N$. We will approximate the inverse Fourier transform by the inverse Discrete Fourier transform (DFT). In order to use Fast Fourier Transform (FFT), we set $\Delta \xi = 1/4$. We choose $\xi_k=-N/8+(k-1)/4$, $k=1,\cdots,N$.

For each $t_m$, we solve the integral equations for $\tilde{j}_n(\xi_k,t_m)$ for all $\xi_l$. We apply the inverse FFT to obtain $j_n(x_l,t_m)$. In computing the convolutions in $\xi$ needed for $L_n$ in \eqref{int-eqn}, we leverage the convolution-multiplication property of Fourier transforms. We compute the multiplication in the spatial domain, then use the FFT to transform the result into the Fourier domain.
Our implementation uses the MATLAB \cite{MATLAB} FFT and IFFT functions.

\subsection{Approximation Errors}
\label{discrete-error}
There are several sources of approximation errors: (i) discretization of the integral equation, (ii) truncating the integration domain in inverse Fourier transform, (iii) evaluating the truncated inverse transform by inverse Discrete Fourier transform (DFT). We will focus on $j_1(x,t)$ and discuss each error source in turn below.

In solving the integral equation, we first need to approximate the kernel in \eqref{kern-def}. Let 
\begin{equation*}
    \overline{K}(\xi; t,\tau) = \int_{\tau}^t e^{-\beta(t-s)}k_1(\xi,s-\tau) ds .
\end{equation*}
Recall $t_m=m\Delta t$ for $m=0,\cdots,M$. We approximate $\overline{K}(\xi; t_m,t_l)$ for $l\leq m$ using the trapezoidal rule, committing an error of $O(\Delta t^2)$. Next, for a given $\xi$, we replace the integral operators in \eqref{j1-int-eqn} and \eqref{int-eqn} by its discrete counterpart
\begin{equation*}
    (K\tilde{j}_1)(\xi,t_m) \approx \sum_{l=0}^m w_l \overline{K}(\xi;t_m,t_l) \tilde{j}_1(\xi,t_l) \Delta t,
\end{equation*}
where $w_0=w_m=1/2$ and $w_l=1$ otherwise. This approximation has an error of $O(\Delta t^2)$. We replace the integral equation \eqref{j1-int-eqn} and \eqref{int-eqn} with
\[
\tilde{j}_{1a_1}(\xi,t_m) + \sum_{l=0}^m w_l \overline{K}(\xi;t_m,t_l) \tilde{j}_{1a_1}(\xi,t_l) \Delta t = L_{n-1}(\xi,t_m).
\]
Following \cite{brunner}, the $l_\infty$ error for $0<m\leq M$ in the approximate solution obtained from the linear system is $O(\Delta t^2)$.

Next, we take the inverse FFT of $\tilde{j}_{1a_1}(\xi_k,t_m)$ to obtain an approximation of $j_1(x_l,t_m)$
\begin{equation*}
    j_{1a_2}(x_l,t_m)=\sum_{k=1}^N \tilde{j}_{1a_1}(\xi_k,t_m) e^{i\xi_k x_l} \Delta \xi.
\end{equation*}
We estimate the total error by writing
\begin{align}
j_1(x_l,t_m) - j_{1a_2}(x_l,t_m) & =
\int_{-\infty}^\infty \tilde{j}_1(\xi,t_m)e^{i\xi x_l}d\xi - \sum_{k=1}^N \tilde{j}_{1a_1}(\xi_k,t_m) e^{i\xi_k x_l} \Delta \xi \nonumber \\
& = \int_{-N/8}^{N/8} \left[\tilde{j}_1(\xi,t_m)
- \tilde{j}_{1a_1}(\xi,t_m) \right] e^{i\xi x_l} d\xi \nonumber\\
&+ \int_{-N/8}^{N/8} \tilde{j}_{1a_1}(\xi,t_m) e^{i\xi x_l} d\xi 
- \sum_{k=1}^N \tilde{j}_{1a_1}(\xi_k,t_m) e^{i\xi_k x_l} \Delta \xi \nonumber \\
&+ \int_{|\xi|>N/8} \tilde{j}_1(\xi,t_m) e^{i\xi x_l} d\xi.
\end{align}
Here, we are assuming there is no aliasing error. It can be seen that there are three error terms.  The first term is $O(\Delta t^2 N)$ since $\tilde{j}_1(\xi,t_m)-\tilde{j}_{1a_1}(\xi,t_m)$ is $O(\Delta t^2)$ irrespective of $\xi$. The second error term can be thought of as the quadrature error and is $O(\Delta\xi^2)$. The final error requires some estimation.

Recall that $g_0(x,t)=d(x,t)E_{x0}(x,0,t)$ and we have assumed that $d(x,t)$ is in $L^2((0,t),H^1(\mathbb{R}))$. From \eqref{boundj1-H1}, we see that $j_1(\cdot,t)\in H^1(\mathbb{R})$. Let $E_3(x,t_m) = \int_{|\xi|>N/2} \tilde{j}_1(\xi,t_m)e^{i\xi x} d\xi$. This implies that
\begin{equation*}
   | E_3(x,t_m) | \leq \int_{|\xi|>N/8} |\tilde{j}_1(\xi,t_m)| d\xi = \int_{|\xi|>N/8} \frac{1}{|\xi|} |\xi \tilde{j}_1(\xi,t_m)| d\xi.
\end{equation*}
We apply Cauchy-Schwartz to the integral and obtain
\begin{align*}
|E_3(x,t_m)| &\leq \left[ \int_{|\xi|>N/8} \frac{1}{\xi^2} d\xi \right]^{1/2}
\left[ \int_{|\xi|>N/8} \xi^2 \tilde{j}_1(\xi,t_m)^2  d\xi\right]^{1/2} \\
& \leq \left[\frac{16}{N} \right]^{1/2} \| j_1(\cdot,t_m) \|_{H^1(\mathbb{R})}^{1/2} .
\end{align*}
Therefore, we can conclude that $|E_3(x_l,t_m)| \leq C/\sqrt{N}$. This, together with the first two error terms, gives a point-wise (in space and time) error estimate.

\subsection{Computational Complexity}
\label{comp-complex}
We consider the computational complexity of evaluating $j_1(x_l,t_m)$ for $l=1,\cdots,N$ and $m=1,\cdots,M$. The cost of evaluating $k_1(\xi,t_m)$ for a fixed $\xi$ using trapezoidal rule is $O(M)$. The discretized kernel in \eqref{kern-def} for $t=t_m$ has $O(M^2)$ entries. Evaluating each entry is $O(M)$ since there is an inner integration. So, the total work of evaluating the discretized operator $K$ is $O(M^3)$. Inverting \eqref{int-eqn} is a matrix inversion for a matrix of size $M$, so the work is $O(M^3)$. This is done for $N$ values of $\xi$. Therefore, to obtain $\tilde{j}_1(\xi_k,t_m)$, the cost is $O(NM^3)$. Finally, to invert the Fourier transform, we use the FFT, which costs $O(N\log N)$. We do this $M$ times for each $t_m$. Therefore, the Fourier inversion cost is $O(MN\log N)$. The total computational cost is $O(NM^3)+O(MN\log N)$. The first term is likely to dominate the second in practice.


\section{Numerical Examples}
\label{num-examples}
The method described in Section 5 has been implemented in MATLAB \cite{MATLAB}. The examples are similar to the ones in \cite{santosa-shi} for comparison. By assumption, the Drude weight perturbation is in $L^2((0,T),H^1(\mathbb{R}))$. We choose $d(x,t)$ to be of the form
\begin{equation*}
    d(x,t) = \frac{D_0}{2}\big( \tanh\left(r(x+a/2)\right)-\tanh\left(r(x-a/2)\right) \big)\cos(\xi_1 x - \omega_1 t),
\end{equation*}
where the $\tanh$ terms form a smoothed rectangle function with width $a$ and sharpness $r$. We set $\Delta t = 0.01575$ and compute from $t=0$ to $t=10$. We also set the following constants: $r = 10$, $a= 4$, $\mu = \epsilon = 1$, $\eta = \sqrt{\mu/\epsilon}$, $D_0 = 0.675$, and $\alpha = 0.05$. 

The approximation is computed using $N=1024$ evenly spaced points from $[-4\pi, 4\pi]$. The plots display the spatial interval of significance, $[-\pi,\pi]$. We simulate for several $\omega_1$ and $\xi_1$ values. Since the background field $j_0(x,t)$ and the perturbation $d(x,t)$ are real. The resulting perturbation solutions $\alpha j_1(x,t)$ and $\alpha^2 j_2(x,t)$ are also real.

\subsection{Propagating current density perturbation}
Let $\xi_1 = \xi_0 = 4$ and $\omega_0 = \omega_1$. We start by considering $j \approx j_0 + \alpha j_1$. 
\begin{figure}[H]
    \centering
\includegraphics[width =0.5\linewidth]{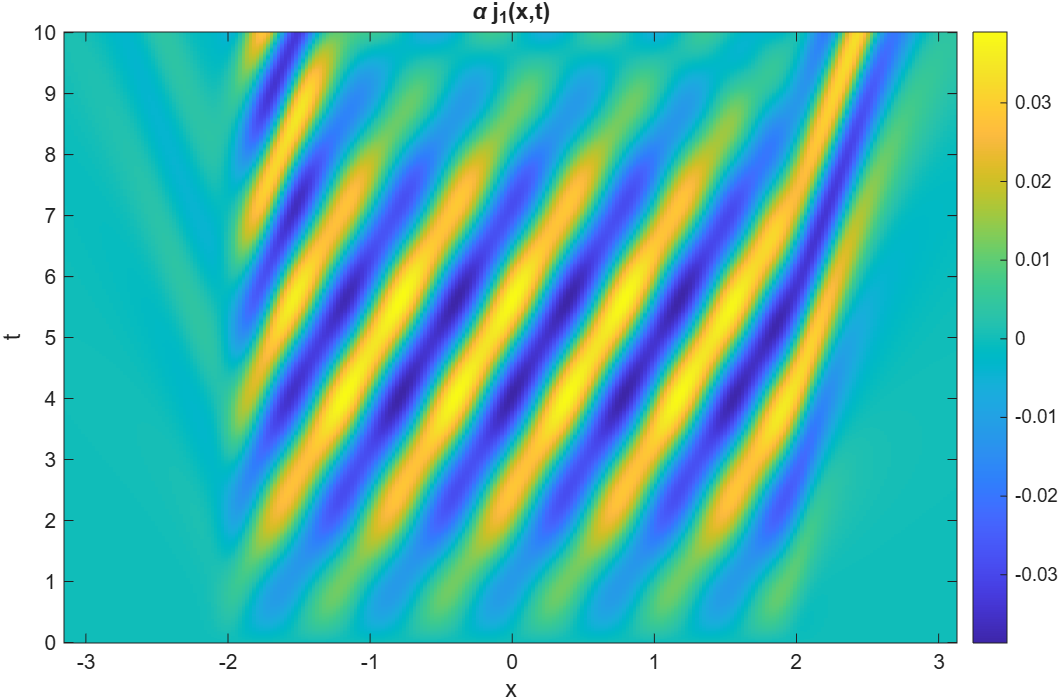}\includegraphics[width=0.5\linewidth]{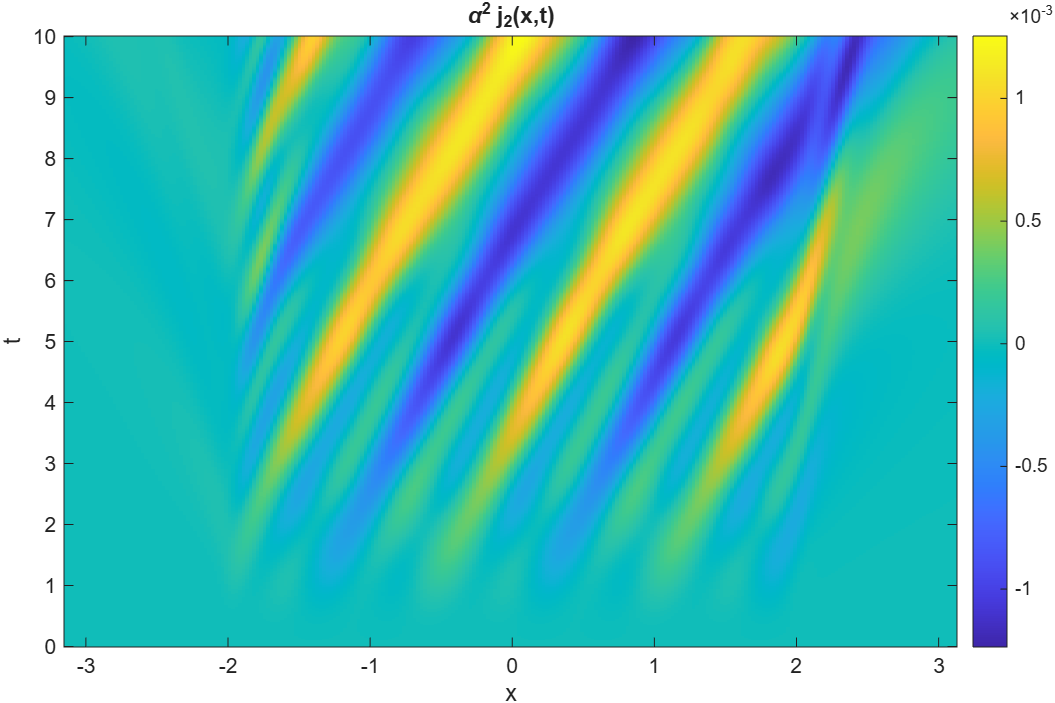}
    \caption{Left: The first order correction term $\alpha j_1(x,t)$ with $\xi_0 = \xi_1 = 4$ and $\omega_0 = \omega_1$. Right: The second order correction term $\alpha^2 j_2(x,t)$.}
    \label{first-order-om1om0}
\end{figure}
We see in Figure \ref{first-order-om1om0} the result is a traveling wave and the correction term is small. Next, add the next $O(\alpha^2)$ term and consider $j \approx j_0 + \alpha j_1 + \alpha^2 j_2$. We expect the contributions from this correction term to also be small. A plot of the second order correction term is displayed. As expected, the correction terms are small. However, it should be noted as time gets larger, the contribution from the correction term grows. 

\subsection{Standing current density perturbations}
Now consider $\omega_0=-\omega_1$. In this case we see standing waves arise in the correction terms. The correction is still concentrated in the support window of $d(x,t)$, as expected, and we see the correction is small.
\begin{figure}[H]
    \centering
    \includegraphics[width=0.5\linewidth]{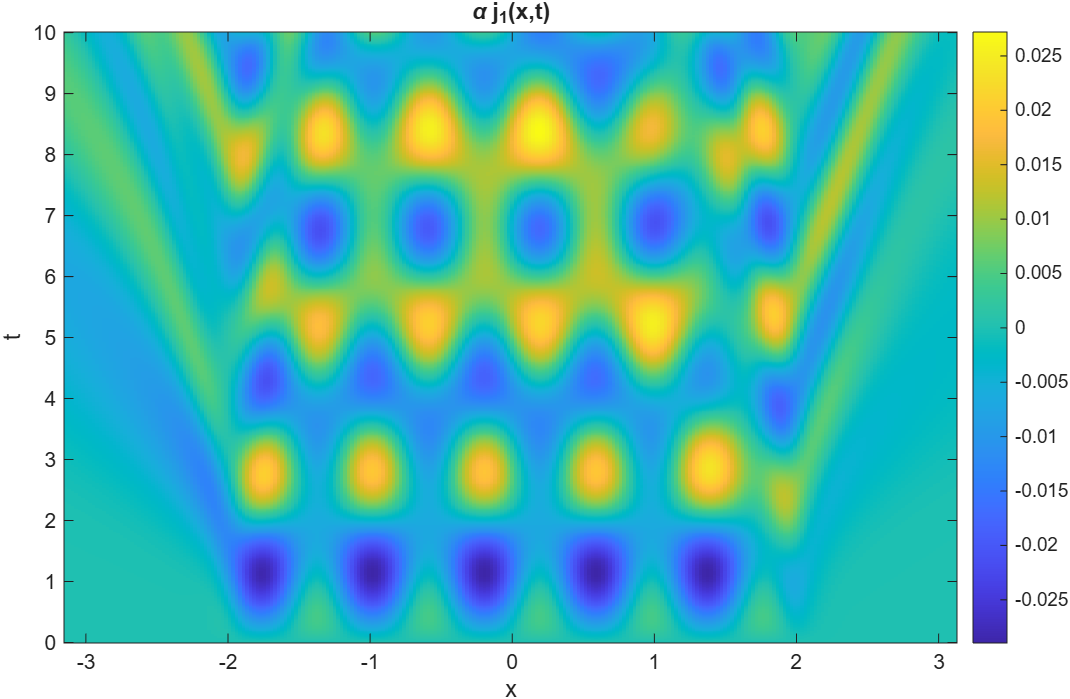}\includegraphics[width=0.5\linewidth]{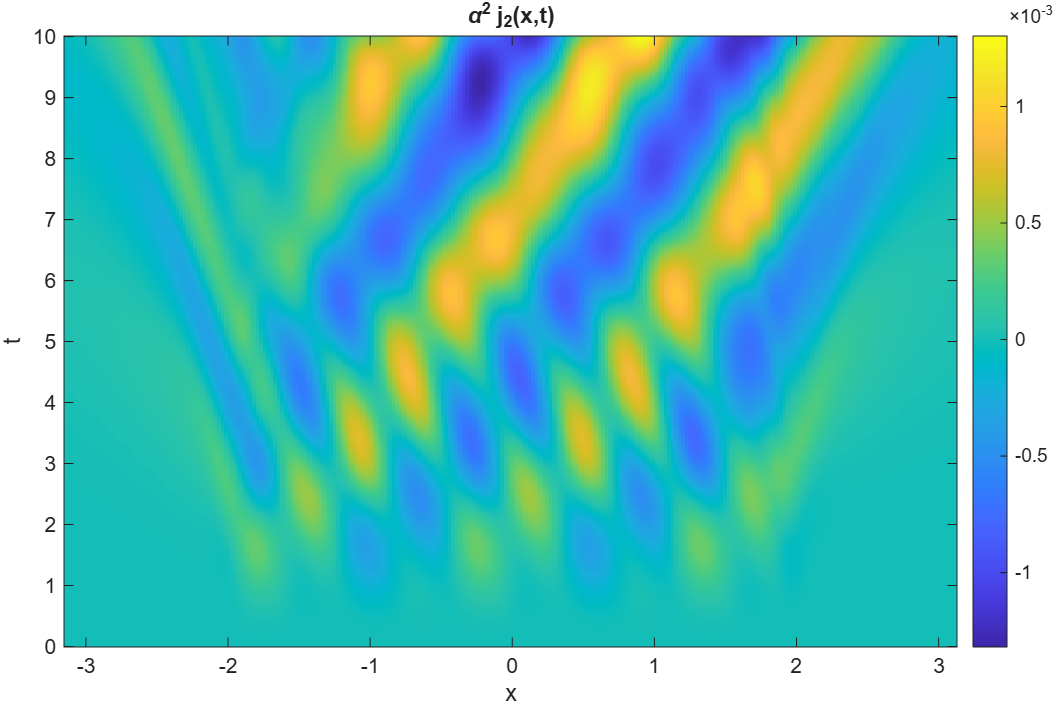}
    \caption{Left: The first order correction term $\alpha j_1(x,t)$ with $\xi_0 = \xi_1 = 4$ and $\omega_0 = -\omega_1$. Right: The second order correction term $\alpha^2 j_2(x,t)$.}
    \label{first-order-om1-om0}
\end{figure}

\subsubsection{Growing current density perturbation}
Consider $\xi_0 = \xi_1 = 4$ and let $\omega_1 = \sqrt{D_0 \xi_0}-\omega_0$. This leads to the behavior of a growing wave in the correction terms. 
\begin{figure}[H]
    \centering
    \includegraphics[width=0.5\linewidth]{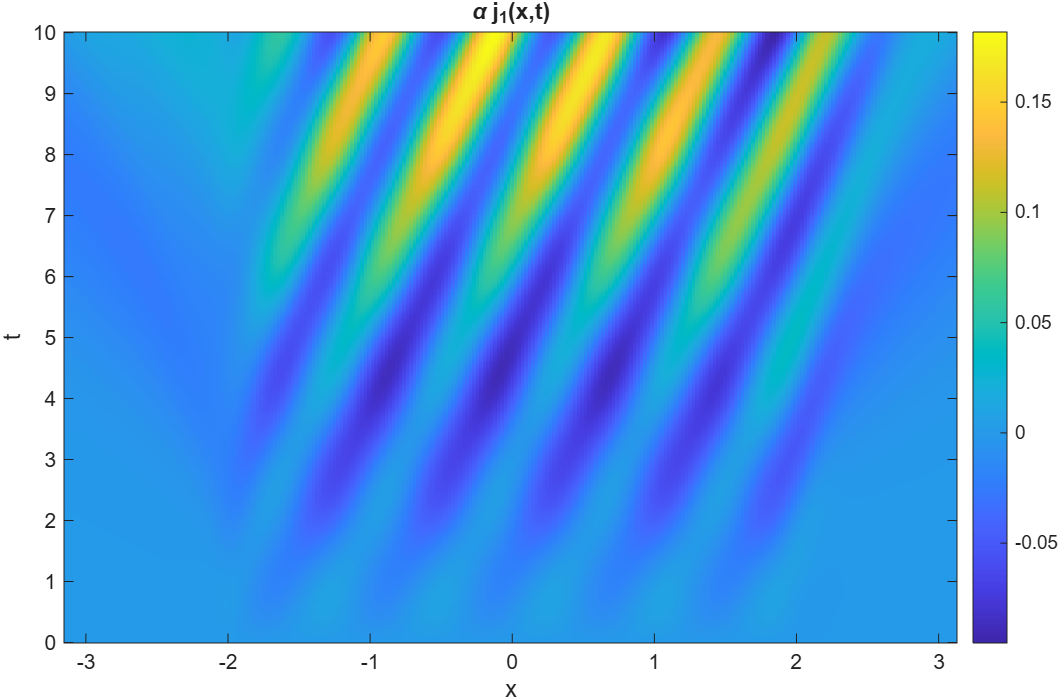}\includegraphics[width=0.5\linewidth]{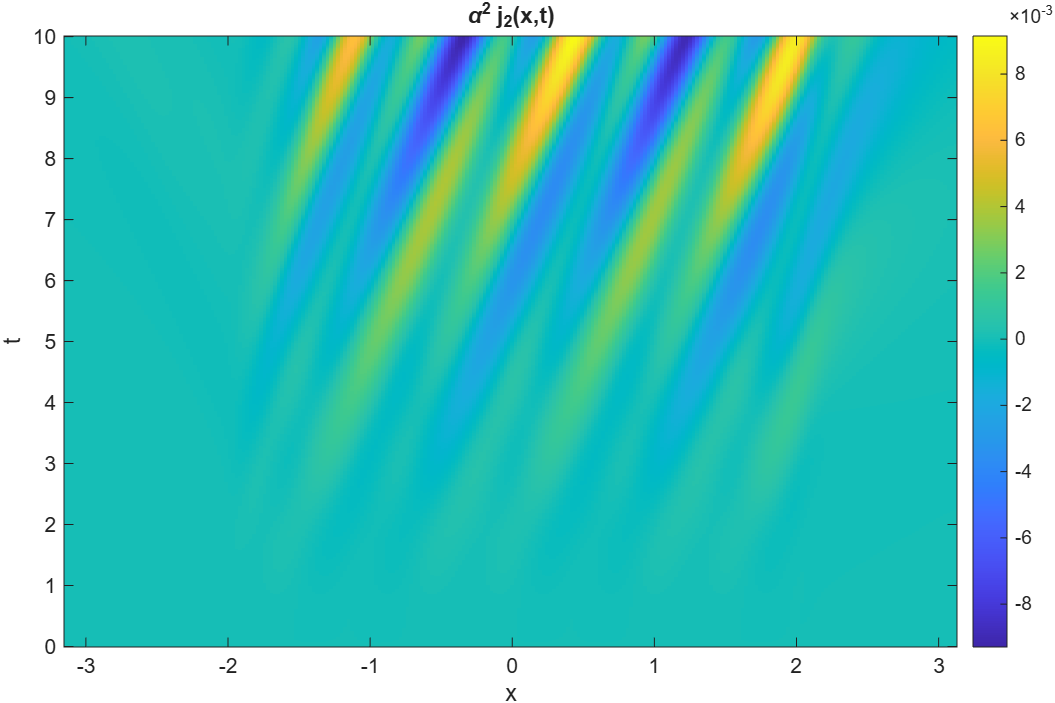}
    \caption{Left: The first order correction term $\alpha j_1(x,t)$ with $\xi_0=\xi_1$ and $\omega_1=\sqrt{D_0\xi_0} - \omega_0$. Right: The second order correction term $\alpha^2 j_2(x,t)$.  }
    \label{growing-wave}
\end{figure}

\section{Solution for Special Form of Drude Weight.}
\label{special-form}
In this section, we consider Drude weight of the form $D(x,t) = D_0 + \alpha d(x,t)$, where
\begin{equation*}
    d(x,t) = D_0\cos{(\xi_1 x-\omega_1 t + \varphi)}.
\end{equation*}
Note that the perturbation in the Drude weight is similar to the one considered in Section 5. However, it should be pointed out that $d(\cdot,t)$ is not in $L^2(\mathbb{R})$. We will present a method to solve for the first order perturbation approximation for $j_1(x,t)$. 

The Fourier transform $j_1(x,t)$ satisfies \eqref{first-ord-eqn} which we rewrite below as
\begin{equation*}
    \frac{d \tilde{j}_1}{dt} + \frac{\eta D_0}{2}\tilde{j}_1 + \frac{\eta D_0}{2}\int_0^t k_1(\xi,t-\tau)\tilde{j_1}(\xi,\tau)d\tau = \tilde{d}(\xi,t)*_\xi \Tilde{E}_{x_0}(\xi,0,t).
\end{equation*}
Inserting the Fourier transforms of $d(x,t)$ and $E_{x_0}(x,0,t)$ in the convolution, we obtain
\begin{multline}
    \frac{d \tilde{j}_1}{dt} + \frac{\eta D_0}{2} \tilde{j}_1 + \frac{\eta D_0}{2} \int_0^t k_1(\xi,t-\tau) \tilde{j}_1 (\xi,\tau) d\tau  =\\ \frac{D_0\gamma_0 \pi}{i\epsilon \omega_0}e^{-i\omega_0 t}\Bigg\{e^{i(\varphi-\omega_1 t)}\delta(\xi-\xi_1-\xi_0)+e^{-i(\varphi-\omega_1 t)}\delta(\xi+\xi_1-\xi_0)\Bigg\}.
    \label{alt_main_eqn}
\end{multline}
The right-hand side is zero except at $\xi$ equal to $\nu_+ = \xi_0+\xi_1$ and $\nu_- = \xi_0-\xi_1$.
We solve \eqref{alt_main_eqn} analytically using the Laplace transform. 

The Laplace transform of a time-dependent function $f(t)$ is indicated by a $\hat{f}(s)$ where
\[
\hat{f}(s) = \int_0^\infty f(t) e^{-st} dt.
\]
Note the Laplace transform of $k_1(\xi,t)$ is given by \cite{wilson-santosa-martin}
\[
\hat{k}_1(\xi,s) = \frac{1}{s}\left[\sqrt{s^2+c^2\xi^2}-s\right].
\]
Then the Laplace transform of \eqref{alt_main_eqn} is:
\begin{multline*}
    s\hat{\tilde{j}}_1 - \tilde{j}_1(\xi,0) + \frac{\eta D_0}{2}\hat{\tilde{j}}_1 + \frac{\eta D_0}{2}\Big(\frac{1}{s}\sqrt{s^2 + c^2\xi^2}-1\Big)\hat{\tilde{j}}_1\\
    = \frac{D_0\gamma_0\pi}{i\epsilon\omega_0}\Bigg\{e^{i\varphi}\frac{\delta(\xi-\nu_+)}{s+i(\omega_0+\omega_1)} + e^{-i\varphi}\frac{\delta(\xi-\nu_-)}{s+i(\omega_0-\omega_1)}\Bigg\} .
\end{multline*}
Let $\beta = \frac{\eta D_0}{2}$. Recall the initial condition gives $\tilde{j}_1(\xi,0) = 0$. Simplifying above to get $\hat{\tilde{j}}_1$ yields
\begin{equation*}
     \hat{\tilde{j}}_1(s,\xi) = \frac{D_0\gamma_0\pi}{i\epsilon\omega_0}\Big(e^{i\varphi}\hat{\tilde{j}}_a(s)\delta(\xi-\nu_+)+e^{-i\varphi}\hat{\tilde{j}}_b(s)\delta(\xi-\nu_-)\Big),
\end{equation*}
where
%
\begin{align}
     \hat{\tilde{j}}_a(s) &= \frac{s^3-\beta s\sqrt{s^2+c^2\nu_+^2}}{(s+i(\omega_0+\omega_1))(s^2-A_+)(s^2-A_-)}\label{ja-eqn}\\
     \hat{\tilde{j}}_b(s) &= \frac{s^3-\beta s\sqrt{s^2+c^2\nu_-^2}}{(s+i(\omega_0-\omega_1))(s^2-B_+)(s^2-B_-)}\label{jb-eqn}.
\end{align}
 Here
\begin{align*}
    A_\pm &= \frac{1}{2}\Big(\beta^2 \pm \sqrt{\beta^4+4\beta^2c^2\nu_+^2}\Big), \\
    B_\pm &= \frac{1}{2}\Big(\beta^2 \pm \sqrt{\beta^4+4\beta^2c^2\nu_-^2}\Big).
\end{align*}
 %
The next task is the Laplace transform inversion of \eqref{ja-eqn} and \eqref{jb-eqn}, which we describe next. 

\subsection{Inverting $\hat{\tilde{j}}_a(s)$}
\label{ja-inversion}
\begin{figure}[t]
\centering
\begin{subfigure}[b]{0.45\textwidth}
         \centering
\begin{tikzpicture}[scale=0.7]
\draw (-3.8,0) -- (4,0) node[right] {\small $\mathrm{Re}(s)$};
\draw (0,-4) -- (0,4) node[above] {\small $\mathrm{Im}(s)$}; 
\draw[thick] (2.7,-1.5) -- (2.7,1.5);
\draw (0,0) -- (2.7,1.5);
\draw[thick] (2.7,1.5) arc [start angle=29, end angle=331, radius=3.09];
\draw[->] (2.7,1.5) arc [start angle=29, end angle=60, radius=3.09];
\draw[->] (2.7,1.5) arc [start angle=29, end angle=145, radius=3.09];
\draw[->] (2.7,1.5) arc [start angle=29, end angle=225, radius=3.09];
\draw[->] (2.7,1.5) arc [start angle=29, end angle=310, radius=3.09];
\node at (-2.4,2.6) {{\small $\Gamma$}};
\node at (1.3,1.1) {{\small $R$}};
\node at (2.9,-0.3) {\small $\gamma$};
\end{tikzpicture}
\caption{The contour for the calculation of $\tilde{ j}_{a_1}$.}\label{bromwich-simple}
\end{subfigure}
     \hfill
\begin{subfigure}[b]{0.45\textwidth}
         \centering   
\begin{tikzpicture}[scale=0.646]
\draw (-4,0) -- (4.4,0) node[right]{\small $\mathrm{Re}(s)$};
\draw (0,-4.35) -- (0,4.4) node[above]{\small $\mathrm{Im}(s)$};
\draw[thick] (3,-1.5) -- (3,1.5);
\draw[->] (3,-1.5) -- ( 3,0.5);
\draw[thick] (3,1.5) arc [start angle=26.57, end angle=88, radius=3.35];
\draw[->] (3,1.5) arc [start angle = 26.57, end angle=57, radius=3.35];
\draw[thick] (3,-1.5) arc [start angle=-26.57, end angle=-88,radius=3.35];
\draw[->] (0.116,-3.34) arc [start angle = -88, end angle = -57, radius=3.35];
\draw [thick] (-0.116,3.34) arc [start angle = 91, end angle=269,radius = 3.35];
\draw [->] (-0.116,3.34) arc [start angle = 91, end angle = 135,radius = 3.35];
\draw [->] (-0.116,3.34) arc [start angle = 91, end angle = 235, radius = 3.35];
\draw[thick] (0.116,2.5) -- (0.116,3.34);
\draw[->] (0.116,3.34) -- (0.116,2.9);
\draw[->] (-0.116,2.5) -- (-0.116,2.9);
\draw[thick] (-0.116,2.5) -- (-0.116,3.34);
\draw[thick] (0.116,2.5) arc[start angle = 66.75, end angle = -246.75,radius=0.275];
\draw[->] (0.116,2.5) arc [start angle = 66.75, end angle = -135,radius = 0.275];
\filldraw[black] (3,1.5) circle (2pt) node[anchor=west]{\small $B$};
\draw (0,0) -- (2.9,1.48);
\node at (1.5,1.2) {{\small $R$}};
\filldraw[black] (3,-1.5) circle (2pt) node[anchor=west]{\small$A$};

\filldraw[black] (0.116,3.34) circle (2pt);
\node at (0.3,3.7) {{\small $C$}};
\filldraw[black] (-0.116,3.34) circle (2pt);
\node at (-0.3,3.7) {\small $E$};
\filldraw[black] (0,2.25) circle (2pt);
\node at (0.2,1.6) {\small $D$};
\draw[thick] (0.116,-2.5) -- (0.116,-3.34);
\draw[->] (0.116,-2.5) -- (0.116,-2.9);
\draw[thick] (-0.116,-2.5) -- (-0.116,-3.34);
\draw[->] (-0.116,-3.34) -- (-0.116,-2.9);
\draw[thick] (0.116,-2.5) arc[start angle = -66.75, end angle = 246.75,radius=0.275];
\draw [->] (-0.116,-2.5) arc [start angle = 246.75, end angle = 135, radius = 0.275];
\filldraw[black] (0.116,-3.34) circle (2pt);
\node at (0.3,-3.7) {\small $H$};
\filldraw[black] (-0.116,-3.34) circle (2pt);
\node at (-0.3,-3.7) {\small $F$};
\filldraw[black] (0,-2.25) circle (2pt);
\node at (0.2,-1.6) {\small $G$};
\node at (3.2,-0.25) {\small $\gamma$};
\end{tikzpicture}
\caption{The contour for calculation of $\tilde{j}_{a_2}$.}\label{bromwich-ja2}
\end{subfigure}
\caption{The contours used in computing $\tilde{j}_{a_1}$ and $\tilde{j}_{a_2}$. In (a), the arc portion is denoted by $\Gamma$. The poles in the expressions are denoted by $p_l$, for $l=1,2,\cdots,5$. The vertical portion must be chosen so that all the poles are to the left of $\gamma$. To obtain the inverse Laplace transform, we will set $R\rightarrow\infty$.}
\end{figure}
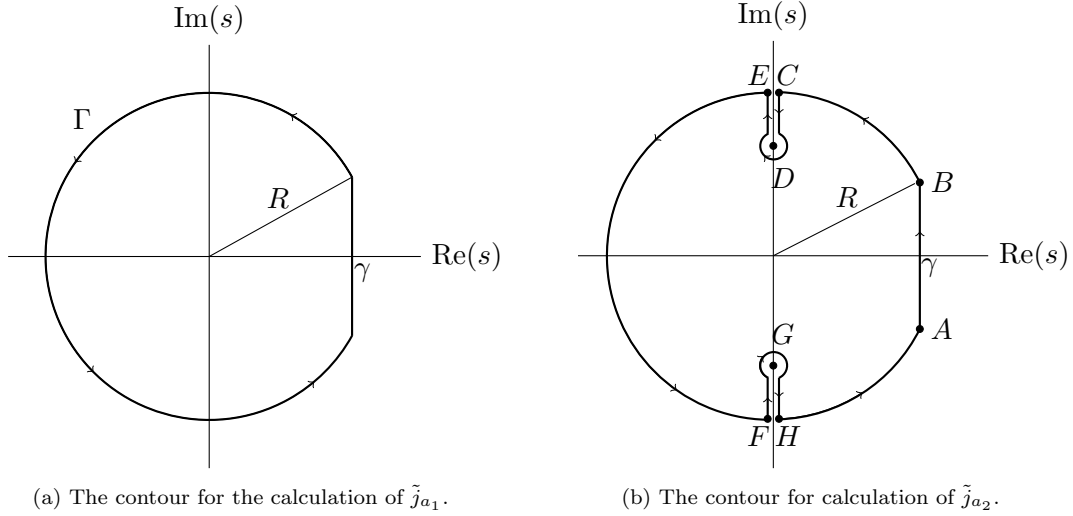

We split $\hat{\tilde{j}}_a(s)$ into two parts
 \begin{align*}
     \hat{\tilde{j}}_{a_1}(s) &= \frac{s^3}{(s+i(\omega_0+\omega_1))(s^2-A_+)(s^2-A_-)},\\
     \hat{\tilde{j}}_{a_2}(s) &= \frac{s\sqrt{s^2+c^2\nu_+^2}}{(s+i(\omega_0+\omega_1))(s^2-A_+)(s^2-A_-)},
 \end{align*}
such that $\hat{\tilde{j}}_a(s) = \hat{\tilde{j}}_{a_1}(s) - \beta\hat{\tilde{j}}_{a_2}(s)$. The term $\hat{\tilde{j}}_{a_1}(s)$ is a rational of two analytic polynomials, there are no branch points in the integral. There are five simple poles: $p_1 = -i(\omega_0+\omega_1)$, $p_2 = -\sqrt{A_+}$, $p_3 = -\sqrt{A_-}$, $p_4 = \sqrt{A_+}$, $p_5 = \sqrt{A_-}$. In order to invert $\hat{\tilde{j}}_{a_1}$ and take into account the poles, we use a Bromwich contour, as shown in Figure \ref{bromwich-simple}, where $\gamma \in \mathbb{R}$ such that $\gamma > \mathrm{Re}(p_l)$ for all $l = 1,2,\cdots,5$. 
Next, apply Cauchy's Integral Theorem and split the contour integral into
\begin{equation*}
    \int_C \hat{\tilde{j}}_{a_1}(s)e^{st}ds = \int_{\gamma -iR}^{\gamma+iR}\hat{\tilde{j}}_{a_1}(s)e^{st}ds + \int_\Gamma \hat{\tilde{j}}_{a_1}(s)e^{st}ds ,
\end{equation*}
where $\Gamma$ is the partial circle of radius, $R$, centered at the origin connecting $\gamma + iR$ and $\gamma-iR$. Using Cauchy's integral formula we have for any closed contour, $C$, 
\begin{equation*}
    \int_C \hat{\tilde{j}}_{a_1}(s)e^{st}ds = 2\pi i \sum_{l=1}^{5} \mathrm{Res}[ \hat{\tilde{j}}_{a_1}(s)e^{st},p_l] .
\end{equation*}
where $p_l$ represents the $l$th pole. As $R\rightarrow \infty$, we can use standard arguments to show that the integral over $\Gamma$ is zero. We are then left with
\[
\tilde{j}_{a_1}(t) = \sum_{l=1}^{5} \mathrm{Res}[ \hat{\tilde{j}}_{a_1}(s)e^{st},p_l] .
\]

Next, consider $\hat{\tilde{j}}_{a_2}(s)$. We note the solution for $\hat{\tilde{j}}_{a_2}(s)$ is not analytic due to the numerator $\sqrt{s^2 + c^2\nu_+^2}$ being a multi-valued function. Therefore, there are two branch points we must consider: $s =\pm ic\nu_+$. Also note we have the same five simple poles as in the $\hat{\tilde{j}}_{a_1}(s)$ case. Due to the combination of poles and branch points we choose the contour shown in Figure \ref{bromwich-ja2}, with the branch cuts to be $[ic\nu_+,iR)$ and $(-iR,-ic\nu_+]$. 

Recall from \cite{wilson-santosa-martin} that $\omega_0 = O(\sqrt{\xi_0})$. We assume that the parameters in the Drude weight perturbation, $\omega_1$ and $\xi_1$, are about the same size as $\omega_0$ and $\xi_0$, respectively. Therefore, we expect $c\nu_+ > (\omega_0+\omega_1)$ and all five poles now lie inside the contour.

Notice this is similar to the Bromwich contour in Figure \ref{bromwich-simple}, but with two key-hole slits. The contributions along the arcs can be shown to vanish as $R \rightarrow \infty$ using the standard argument. One can also argue that the contributions around the points $D$ and $G$ are zero.

The computation of the contributions due to the vertical branch cuts is nontrivial, and we provide some details here. First, consider the integral over the segment $CD$. We have
\begin{equation*}
    I_{CD} =
    \int_{iR}^{ic\nu}\frac{s(s-ic\nu_+)^{1/2}(s+ic\nu_+)^{1/2}e^{st}}{(s+i(\omega_0+\omega_1))(s-\sqrt{A_+})(s-\sqrt{A_-})(s+\sqrt{A_+})(s+\sqrt{A_-)}}ds.
\end{equation*}
The segment $CD$ is on the positive real half-space. We have eight terms in the integrand and write each term in the form of $f(s)=|f(s)|e^{i\mathrm{arg}(f(s))}$. Since there is a branch cut, we need to carefully restrict the values of the argument to avoid being on a different branch. Let $\mathrm{arg}(s-ic\nu_+) \in (-3\pi/2,\pi/2]$ and let $\mathrm{arg}(s+ic\nu_+) \in (-\pi/2,3\pi/2]$. The arguments, restricted to $(-\pi,\pi]$, are shown in Table \ref{args-ja}.
\begin{table}[t]
\begin{center}
\begin{tabular}{ |p{2.6cm}|p{2cm}|p{2cm}| }
 \hline
 Term & argument for $CD$ & argument for $DE$\\
 \hline
 $s-ic\nu$& $\pi/2$ & $-3\pi/2$\\
 $s+ic\nu$& $\pi/2$ & $\pi/2$\\
 $s$ & $\pi/2$ & $\pi/2$\\
 $s+i(\omega_0+\omega_1)$ & $\pi/2$ & $\pi/2$\\
 $s+\sqrt{A_+}$ & $\pi/2$ & $\pi/2$\\
 $s-\sqrt{A_+}$ & $\pi/2$ & $\pi/2$\\
 $s+\sqrt{A_-}$ & $\pi/2$ & $\pi/2$\\
 $s-\sqrt{A_-}$ & $\pi/2$ & $\pi/2$\\
 \hline
\end{tabular}
\end{center}
\caption{Table of arguments for each part of the expressions for $I_{CD}$ and $I_{DE}$.}
\label{args-ja}
\end{table}

Taking into account the argument of each term, we have
\begin{equation*}
 I_{CD} = \int_{iR}^{ic\nu_+}\frac{e^{\frac{-3\pi i}{2}}|s\|s-ic\nu_+|^{1/2}|s+ic\nu_+|^{1/2}e^{st}}{|s+i(\omega_0+\omega_1)\|s-\sqrt{A_+}\|s-\sqrt{A_-}\|s+\sqrt{A_+}\|s+\sqrt{A_-|}}ds.
\end{equation*}
Make the following change of variables $s = iy$ and obtain
\begin{equation*}
    I_{CD} = \int_{c\nu_+}^{R}\frac{|y\|y-c\nu_+|^{1/2}|y+c\nu_+|^{1/2}e^{iyt}}{|y+(\omega_0+\omega_1)\|iy-\sqrt{A_+}\|iy-\sqrt{A_-}\|iy+\sqrt{A_+}\|iy+\sqrt{A_-|}}dy.\\
\end{equation*}
For $DE$, the only difference is that we accumulate an argument of $-5\pi/2$, which means that $I_{DE}=I_{CD}$.

Similar calculations lead to
\begin{align*}
    I_{FG} &= \int_{-R}^{-c\nu_+}\frac{|y\|y-c\nu_+|^{1/2}|y+c\nu_+|^{1/2}e^{iyt}}{|y+(\omega_0+\omega_1)\|iy-\sqrt{A_+}\|iy-\sqrt{A_-}\|iy+\sqrt{A_+}\|iy+\sqrt{A_-|}}dy \\
& = I_{GH}.
\end{align*}

Combining all of the contributions gives
\begin{equation*}
    \tilde{j}_a(t) = \sum_{l=1}^{5} 
    \text{Res}[\hat{\tilde{j}}_{a}(s)e^{st}, p_{l}]
     - \frac{1}{\pi i}\left( \lim_{R\rightarrow\infty}I_{CD} - \lim_{R\rightarrow\infty}I_{FG} \right).
\end{equation*}
It should be noted that $p_4$ is real positive and appear as exponentially growing terms in the residues. However, when added with $p_2$, the residues cancel each other. Therefore, the residues contribute only oscillatory or decaying terms. The integrals are not available in closed-form and are evaluated numerically.

\subsection{Inverting $\hat{\tilde{j}}_b(s)$}
\label{jb-inversion}
Inverting $\hat{\tilde{j}}_b(s)$ is a bit more involved as there are special cases to consider and there is a possibility that one of the poles lies outside the contour. For details of this computation, we refer the reader to \ref{detail-jb-inversion}. The formulas for $\tilde{j}_b(t)$ is given in \eqref{jb-casei}, \eqref{jb-caseii}, \eqref{jb-caseiii}, \eqref{jb-caseiva}, \eqref{jb-caseivb} and \eqref{jb-caseivc}.

\subsection{Inverting the Fourier Transform}
\label{fourier-inversion}
After inverting the Laplace transform, we invert the Fourier transform. Notice due to the way $\tilde{j}_a(t)$ and $\tilde{j}_b(t)$ are defined, they are not dependent on $\xi$, so $\tilde{j}_a(t) = j_a(t)$ and $\tilde{j}_b(t) = j_b(t)$. Then
\begin{align*}
    \tilde{j}_1(\xi,t) &= \frac{D_0\gamma_0\pi}{i\epsilon\omega_0}\Big(e^{i\varphi}j_a(t)\delta(\xi-\nu)+e^{-i\varphi}j_b(t)\delta(\xi-\zeta)\Big).
\end{align*}
Since the inverse Fourier transform is linear, we only invert the $\delta$-functions. Thus we have
\begin{equation}
    j_1(x,t) = \frac{D_0\gamma_0}{2i\epsilon\omega_0}\Big(e^{i\varphi}j_a(t)e^{i(\xi_0+\xi_1)x}+e^{-i\varphi}j_b(t)e^{i(\xi_0-\xi_1)x}\Big).
\end{equation}

\subsection{Numerical examples}
\label{special-num-exps}
We perform the same numerical experiments as presented in \citep{santosa-shi}. For the following experiments we let $D(x,t) = D_0 + \alpha D_0 \cos(\xi_1x -\omega_1 t)$, $\xi_0 = 4$, and $D_0 = 0.675$. This leads to $\omega_0 = 1.1377$. We show the function $j_1(x,t)$.
\subsubsection{Traveling perturbational current density}
\begin{figure}[H]
    \centering
    \includegraphics[width=0.7\linewidth]{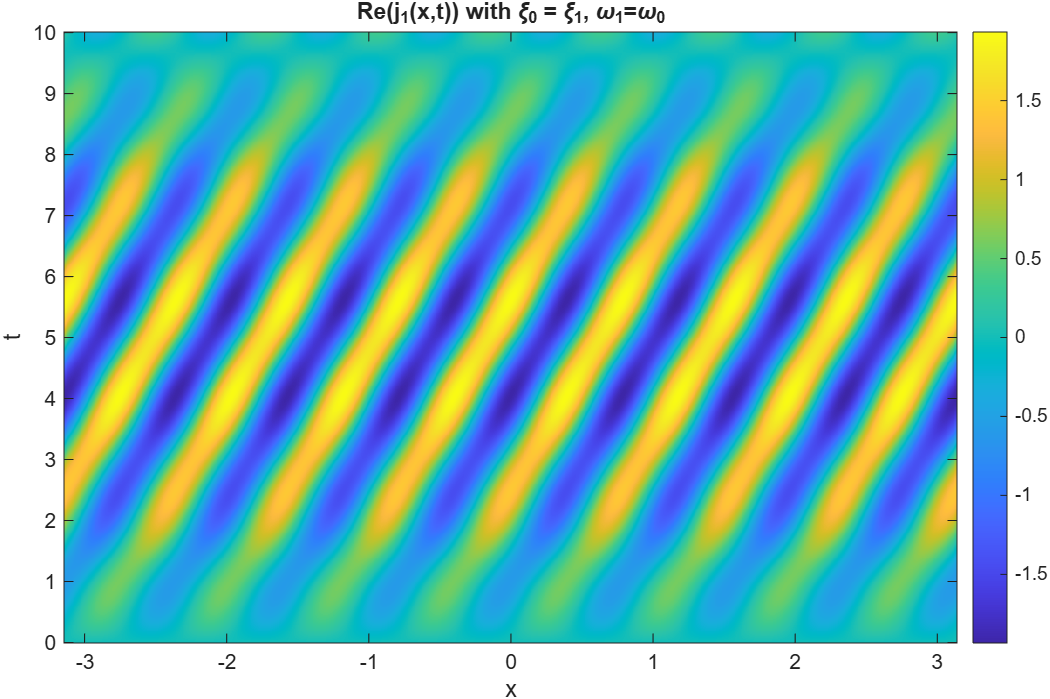}
    
    \caption{Plot of the perturbation term, $j_1$. We notice the perturbation term takes the form of a traveling wave.}
    \label{new-propag}
\end{figure}

Let $\xi_0 = \xi_1$ and let $\omega_0 = \omega_1$. With these parameters, we see in Figure \ref{new-propag} a traveling wave arise in the perturbation correction term. The solution should be compared with Figure 6 in \cite{santosa-shi}.

\begin{figure}[H]
    \centering
    \includegraphics[width=0.7\linewidth]{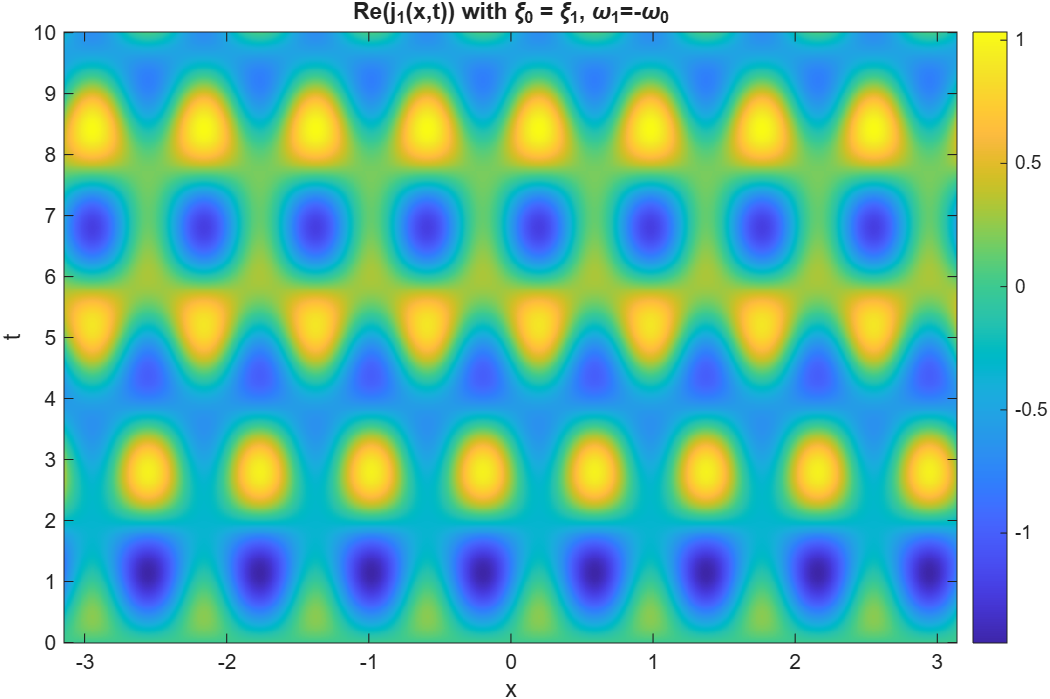}
    
    \caption{Plot of the perturbation term, $j_1$. In this case, the perturbation term takes the form of a standing wave.}
    \label{new-standing}
\end{figure}
\subsubsection{Standing perturbational current density}
Let $\xi_0 = \xi_1$ and let $\omega_0 = -\omega_1$. With these parameters, we see in Figure \ref{new-standing} standing waves arise in the perturbation term. Compare with Figure 7 in \cite{santosa-shi}.

\subsubsection{Growing perturbational current density}
\begin{figure}[H]
    \centering
    \includegraphics[width=0.7\linewidth]{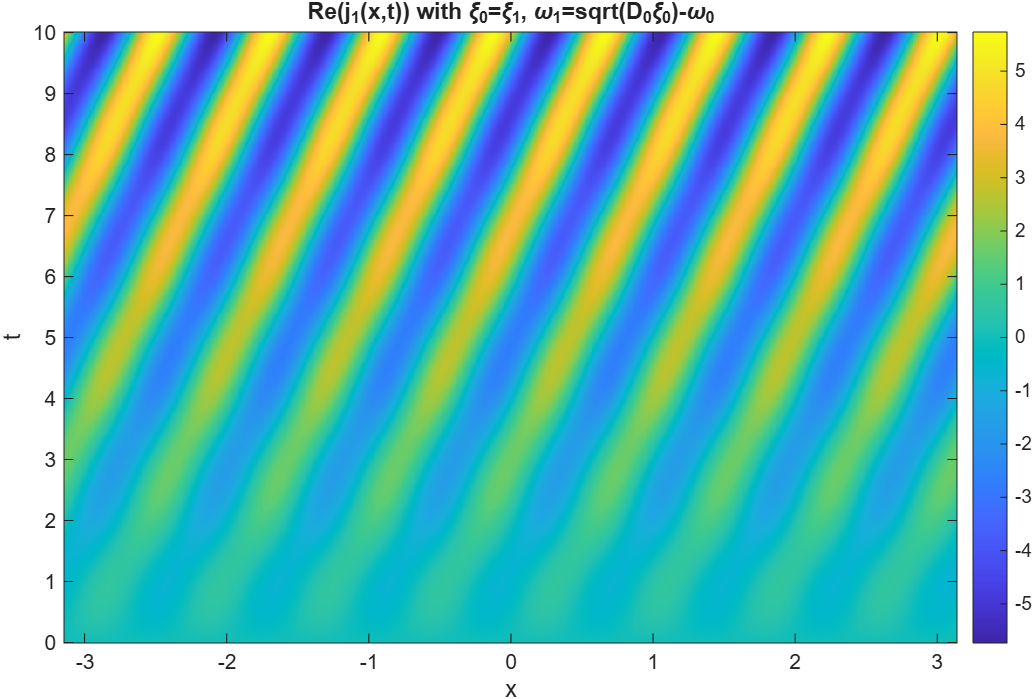}
    \caption{Plot of the perturbation term, $j_1$. We notice the perturbation term takes the form of a growing, traveling wave.}
    \label{new-growing}
\end{figure}

Let $\xi_0=\xi_1$ and $\omega_1 = \sqrt{D_0\xi_0}-\omega_0$. These parameters allow a traveling wave that is growing in amplitude (See Figure \ref{new-growing}) in the perturbation term. Compare with Figure 8 in \cite{santosa-shi}

\section{Discussion}
\label{Discussion}
We propose two methods for efficient simulation of perturbed surface plasmons on graphene where the background constant Drude weight is perturbed in time and space. The first method formulates the governing PIDE into a hierarchy of Volterra integral equations parameterized by wave number $\xi$. We show existence, uniqueness, and spatial regularity of the perturbation. The developed discrete numerical scheme scales with $O(NM^3 + MN \log N)$ complexity where $M$ is the number of time steps and $N$ is the number of spatial samples. The method offers high accuracy for bounded spatiotemporal perturbations and is more efficient than direct discretization which needs to deal with absorbing boundary conditions. The second method is based on inverting the transforms via contour integration and quadrature. It works on a special form of Drude weight perturbation. Numerical examples, similar to those in \cite{santosa-shi}, are provided to demonstrate the methods. These approximate methods offer valuable design and predictive tools for integrated photonic devices that leverages plasmons on 2D materials.

\bibliographystyle{elsarticle-num-names} 
\bibliography{rudge-santosa.bib}

\appendix
\setcounter{figure}{0}
\section{Inversion of $\hat{\tilde{j}}_b(s)$}
\label{detail-jb-inversion}
We start with
\[
\hat{\tilde{j}}_{b}(s) = \frac{s^3-\beta s \sqrt{s^2+c^2\nu_-^2}}{(s+i(\omega_0-\omega_1))(s^2-B_+)(s^2-B_-)}.
\]
Since $\nu_- = \xi_0-\xi_1$, it can take on zero as a value. Therefore, we must consider four cases: (i) $\nu_- = 0$ and $\omega_0-\omega_1 = 0$; (ii) $\nu_- = 0$ and $\omega_0-\omega_1 \neq 0$; (iii) $\nu_- \neq 0$ and $\omega_0-\omega_1 = 0$; and (iv) $\nu_- \neq 0$ and $\omega_0-\omega_1 \neq 0$. 

\medskip
\noindent Case (i): $\nu_- = 0$ and $\omega_0-\omega_1 = 0$.
    
Since $\nu_- = 0$, we have $B_- = 0$. In addition, with $\omega_0-\omega_1=0$, $\hat{\tilde{j}}_b(s)$ simplifies to
    \begin{equation*}
        \hat{\tilde{j}}_{b}(s) = \frac{1}{s(s+\beta)} = \frac{1}{\beta}\left(\frac{1}{s}-\frac{1}{s+\beta}\right),
    \end{equation*}
by partial fractions. This leads to
    \begin{equation}
        \tilde{j}_b(t) = \frac{1}{\beta}\left(1-e^{-\beta t}\right). \label{jb-casei}
    \end{equation}

\medskip
\noindent Case (ii) $\nu_- = 0$ and $\omega_0-\omega_1 \neq 0$.

We again have $B_- = 0$, so $\hat{\tilde{j}}_b(s)$ simplifies to
    \begin{align*}
        \hat{\tilde{j}}_{b}(s) &= \frac{1}{(s+i(\omega_0-\omega_1))(s+\beta)}.
    \end{align*}
The inverse Laplace transform is 
    \begin{equation}
        \tilde{j}_b(t) =  \frac{1}{\beta-i(\omega_0-\omega_1)}\left(e^{-i(\omega_0-\omega_1)t}-e^{-\beta t} \right)
        . \label{jb-caseii}
    \end{equation}

\medskip
    \noindent Case (iii) $\nu_- \neq 0$ and $\omega_0-\omega_1 = 0$.

In this case, we proceed as in the $\hat{\tilde{j}}_a(s)$ case. With $\omega_0-\omega_1 = 0$, we have
\[
\hat{\tilde{j}}_{b}(s) = \frac{s^2-\beta\sqrt{s^2+c^2\nu_-^2}}{s(s^2-B_+)(s^2-B_-)}.
\]
The expression has four simple poles at $q_{1,2} = \pm \sqrt{B_+}$ and $q_{3,4}=\pm \sqrt{B_-}$. Then, as before, we use a Bromwich contour in Figure A.1(a) and Cauchy's integral theorem. Together with the contribution along the branch cuts, we arrive at
\begin{equation}
    \tilde{j}_{b}(t) = \sum_{l=1}^{4}\text{Res}[\hat{\tilde{j}}_{b}(s)e^{st},q_l] -\frac{1}{\pi i}\left(
    \lim_{R\rightarrow\infty} I_{C'D'} + \lim_{R\rightarrow\infty} I_{F'G'} \right), \label{jb-caseiii}
\end{equation}
where
\begin{align*}
I_{C'D'} &= \int_{c\nu_-}^R \frac{\sqrt{y^2-c^2\nu_-^2}e^{iyt}}{(y^2+B_+)(y^2+B_-)}dy,\\
I_{F'G'} &= -\int_{-R}^{-c\nu_-} \frac{\sqrt{y^2-c^2\nu_-^2}e^{iyt}}{(y^2+B_+)(y^2+B_-)}dy.
\end{align*}

\begin{figure}[t]
\centering
\begin{subfigure}[b]{0.45\textwidth}
       \centering
\begin{tikzpicture}[scale=0.658]
    \draw (-4,0) -- (4,0) node[right]{\small $\mathrm{Re}(s)$};
    \draw (0,-4.2) -- (0,4.2) node[above]{\small $\mathrm{Im}(s)$};
    \draw[thick] (3,-1.5) -- (3,1.5);
    \draw[->] (3,-1.5) -- ( 3,0.5);
    \draw[thick] (3,1.5) arc [start angle=26.57, end angle=88, radius=3.35];
    \draw[->] (3,1.5) arc [start angle = 26.57, end angle=57, radius=3.35];
    \draw[thick] (3,-1.5) arc [start angle=-26.57, end angle=-88,radius=3.35];
    \draw[->] (0.116,-3.34) arc [start angle = -88, end angle = -57, radius=3.35];
    \draw [thick] (-0.116,3.34) arc [start angle = 91, end angle=269,radius = 3.35];
    \draw [->] (-0.116,3.34) arc [start angle = 91, end angle = 135,radius = 3.35];
    \draw [->] (-0.116,3.34) arc [start angle = 91, end angle = 235, radius = 3.35];
    \draw[thick] (0.116,1.5) -- (0.116,3.34);
    \draw[->] (0.116,3.34) -- (0.116,2.4);
    \draw[->] (-0.116,1.5) -- (-0.116,2.34);
    \draw[thick] (-0.116,1.5) -- (-0.116,3.34);
    \draw[thick] (0.116,1.5) arc[start angle = 66.75, end angle = -246.75,radius=0.275];
    \draw[->] (0.116,1.5) arc [start angle = 66.75, end angle = -135,radius = 0.275];
    \filldraw[black] (3,1.5) circle (2pt) node[anchor=west]{{\small $B'$}};
    \draw (0,0) -- (3,1.5);
\node at (1.5,1.2) {{\small $R$}};
    \filldraw[black] (3,-1.5) circle (2pt) node[anchor=west]{{\small $A'$}};
    \filldraw[black] (0.116,3.34) circle (2pt);
    \node at (0.34,3.7) {{\small $C'$}};
    \filldraw[black] (-0.116,3.34) circle (2pt);
    \node at (-0.34,3.7) {{\small $E'$}};
    \filldraw[black] (0,1.25) circle (2pt);
    \node at (0.6,1.0) {{\small $D'$}};
    \draw[thick] (0.116,-1.5) -- (0.116,-3.34);
    \draw[->] (0.116,-1.5) -- (0.116,-2.9);
    \draw[thick] (-0.116,-1.5) -- (-0.116,-3.34);
    \draw[->] (-0.116,-3.34) -- (-0.116,-1.9);
    \draw[thick] (0.116,-1.5) arc[start angle = -66.75, end angle = 246.75,radius=0.275];
    \draw [->] (-0.116,-1.5) arc [start angle = 246.75, end angle = 135, radius = 0.275];
    \filldraw[black] (0.116,-3.34) circle (2pt);
    \node at (0.34,-3.7) {{\small $H'$}};
    \filldraw[black] (-0.116,-3.34) circle (2pt);
    \node at (-0.34,-3.7) {{\small $F'$}};
    \filldraw[black] (0,-1.25) circle (2pt);
    \node at (0.6,-1.0) {{\small $G'$}};
    \node at (3.2,-0.3) {$\gamma$};
    \end{tikzpicture}
\caption{The contour for the calculation of $\tilde{ j}_{b}$ when $\omega_0=\omega_1$.}\label{bromwich-jb-caseiii}
\end{subfigure}
 \hfill 
\begin{subfigure}[b]{0.45\textwidth}
         \centering   
\begin{tikzpicture}[scale=0.646,transform canvas={yshift=2.8cm}]
\draw (-4,0) -- (4,0) node[right]{$\mathrm{Re}(s)$};
\draw (0,-4.2) -- (0,4.2) node[above]{$\mathrm{Im}(s)$};
\draw[thick] (3,-1.5) -- (3,1.5);
\draw[->] (3,-1.5) -- ( 3,0.5);
\draw[thick] (3,1.5) arc [start angle=26.57, end angle=88, radius=3.35];
\draw[->] (3,1.5) arc [start angle = 26.57, end angle=57, radius=3.35];
\draw[thick] (3,-1.5) arc [start angle=-26.57, end angle=-88,radius=3.35];
\draw[->] (0.116,-3.34) arc [start angle = -88, end angle = -57, radius=3.35];
\draw [thick] (-0.116,3.34) arc [start angle = 91, end angle=269,radius = 3.35];
\draw [->] (-0.116,3.34) arc [start angle = 91, end angle = 135,radius = 3.35];
\draw [->] (-0.116,3.34) arc [start angle = 91, end angle = 235, radius = 3.35];
\draw[thick] (0.116,1.5) -- (0.116,3.34);
\draw[->] (0.116,3.34) -- (0.116,2.5);
\draw[->] (-0.116,1.5) -- (-0.116,2.5);
\draw[thick] (-0.116,1.5) -- (-0.116,3.34);
\draw[thick] (0.116,1.5) arc[start angle = 66.75, end angle = -246.75,radius=0.275];
\draw[->] (0.116,1.5) arc [start angle = 66.75, end angle = -135,radius = 0.275];
\filldraw[black] (3,1.5) circle (2pt) node[anchor=west]{{\small $B'$}};
\draw (0,0) -- (3,1.5);
\node at (1.5,1.2) {{\small $R$}};
\filldraw[black] (3,-1.5) circle (2pt) node[anchor=west]{{\small $A'$}};
\filldraw[black] (0.12,3.34) circle (2pt);
\node at (0.35,3.7) {{\small $C'$}};
\filldraw[black] (-0.116,3.34) circle (2pt);
\node at (-0.35,3.7) {{\small $E'$}};
\filldraw[black] (0,1.25) circle (2pt);
\node at (0.6,1) {{\small $D'$}};
\draw[thick] (0.116,-1.5) -- (0.116, -2.25);
\draw[thick] (0.116,-2.75) -- (0.116,-3.34);
\draw[->] (0.116,-1.5) -- (0.116,-1.8);
\draw[->] (0.116,-2.75) -- (0.116,-3);
\draw[thick] (-0.116,-1.5) -- (-0.116, -2.25);
\draw[thick] (-0.116,-2.75) -- (-0.116,-3.34);
\draw[->] (-0.116,-3.34) -- (-0.116,-2.9);
\draw[->] (-0.116,-2.25) -- (-0.116,-1.75);
\draw[thick] (0.116,-2.75) arc[start angle = -66.75,end angle = 66.75, radius=0.275];
\draw[->] (0.116,-2.25) arc[start angle = 66.75, end angle = 0, radius=0.275];
\draw[thick] (-0.116,-2.77) arc[start angle = -86.75,end angle = 86.75, radius=0.275];
\draw[->] (-0.116,-2.77) arc[start angle = -86.75,end angle = 0, radius=0.275];
\draw[thick] (0.116,-1.5) arc[start angle = -66.75, end angle = 246.75,radius=0.275];
\draw [->] (-0.116,-1.5) arc [start angle = 246.75, end angle = 135, radius = 0.275];
\filldraw[black] (0.116,-3.34) circle (2pt);
\node at (0.35,-3.7) {{\small $L'$}};
\node at (-0.35,-2.7) {{\scriptsize $G'$}};
\node at (-0.4,-2.2) {{\scriptsize $H'$}};
\node at (0.41,-2.8) {{\scriptsize $K'$}};
\filldraw[black] (0,-2.5) circle (2pt);
\node at (0.44,-2.1) {{\scriptsize $J'$}};
\filldraw[black] (-0.116,-3.34) circle (2pt);
\node at (-0.35,-3.7) {{\small $F'$}};
\filldraw[black] (0,-1.25) circle (2pt);
\node at (0.6,-1) {{\small $I'$}};
\node at (3.2,-0.3) {$\gamma$};
\end{tikzpicture}
\caption{The contour with two key-hole branch cuts and a notch for the general case.}\label{bromwich-jb-caseiv}
\end{subfigure}
\caption{The Bromwich contours for $\hat{\tilde{j}}_{b}$ for Cases (iii) and (iv)}
\end{figure}

\noindent Case (iv) $\nu_- \neq 0$ and $\omega_0-\omega_1 \neq 0$.
\medskip

As we did in Subsection \ref{ja-inversion}, we split $\hat{\tilde{j}}_b(s)$ in \eqref{jb-eqn} into two parts
\begin{align*}
    \hat{\tilde{j}}_{b_1}(s) &= \frac{s^3}{(s+i(\omega_0-\omega_1))(s^2-B_+)(s^2-B_-)},\\
    \hat{\tilde{j}}_{b_2}(s) &= \frac{s\sqrt{s^2+c^2\nu_-^2}}{(s+i(\omega_0-\omega_1))(s^2-B_+)(s^2-B_-)},
\end{align*}
so that $\hat{\tilde{j}}_b(s) = \hat{\tilde{j}}_{b_1}(s)-\beta\hat{\tilde{j}}_{b_2}(s)$.
    
We see, as before, for $\hat{\tilde{j}}_{b_1}(s)$ there are no branch cuts, only five simple poles, so we can solve it by applying Cauchy's integral theorem and take the sum of the residuals. We note that    $\hat{\tilde{j}}_{b_2}(s)$ has two branch points at $\pm ic\nu_-$ and five simple poles. We must separately consider the cases where $c\nu_-$ is greater than, less than, or equal to $\omega_0-\omega_1$.

\medskip
\noindent
Sub-case (a): $c\nu_- > \omega_0 - \omega_1$

The evaluation of this case is similar to the analysis in Case (iii), with the addition of a fifth pole. We use the contour in Figure \ref{bromwich-jb-caseiii}(a). The result is
\begin{equation}
    \tilde{j}_b(t) = \sum_{l=1}^{5} \text{Res}[\hat{\tilde{j}}_b(s)e^{st},q_{l}]  - \frac{1}{\pi i}\left(\lim_{R\rightarrow\infty} I_{F'G'} + \lim_{R\rightarrow\infty} I_{C'D'} \right), 
    \label{jb-caseiva}
\end{equation}
where $q_{l}$ are the $l$th poles of $\hat{\tilde{j}}_{b}(s)$ and 
\begin{align*}
    I_{F'G'} &= -\int_{-R}^{-c\nu_-}\frac{y\sqrt{y^2-c^2\nu_-^2}e^{iyt}}{(y+\omega_0-\omega_1)(y^2+B_+)(y^2+B_-)}dy, \\
    I_{C'D'} &= \int_{c\nu_-}^{R}\frac{y\sqrt{y^2-c^2\nu_-^2}e^{iyt}}{(y+\omega_0-\omega_1)(y^2+B_+)(y^2+B_-)}dy. 
\end{align*}

\medskip
\noindent Sub-case (b): $c\nu_- < \omega_0-\omega_1$

The pole at $-i(\omega_0-\omega_1)$ lies below the branch point $-ic\nu_-$. We use the contour in Figure \ref{bromwich-jb-caseiv} where one can see an additional notch around the pole in the lower branch cut. The contribution of the integration over the arcs can be shown to be zero as $R\rightarrow\infty$.

We start with the calculation of the upper branch cut, i.e. the contributions of $C'D'$ and $D'E'$. This is quite straight-forward
\begin{equation*}
    I_{C'D'} = \int_{c\nu_-}^{R} \frac{y(y^2-c^2\nu_-^2)^{1/2}e^{iyt}}{(y+(\omega_0-\omega_1))(y^2+B_+)(y^2+B_-)}dy  =I_{D'E'} .
\end{equation*}

Next, consider the calculation of the lower branch cut. Since we now have a notch around the simple pole, this cut has more components to consider. Let the notch around the simple pole have radius $\varepsilon$. We consider the components in the limit as $\varepsilon \rightarrow 0$. Let $\omega = \omega_0-\omega_1$. Looking at the contributions going up on the left side, we have
\begin{align*}
    I_{F'G'}&=\int_{-R}^{-\omega-\varepsilon}\frac{(iy)(-y^2+c^2\nu_-^2)^{1/2}e^{iy t}}{(iy+i\omega)(-y^2-B_+)(-y^2-B_-)}idy,\\
    I_{G'H'}&= \int_{-\pi/2}^{\pi/2}\frac{(-i\omega+\varepsilon e^{i\theta})((-i\omega+\varepsilon e^{i\theta})^2+c^2\nu_-^2))^{1/2}e^{(-i\omega+\varepsilon e^{i\theta})t}}
    {((-i\omega+\varepsilon e^{i\theta})^2-B_+)((-i\omega+\varepsilon e^{i\theta})^2-B_-)}id\theta,\\
    I_{H'I'} &= \int_{-\omega+\varepsilon}^{-c\nu_-}\frac{(iy)(-y^2+c^2\nu_-^2)^{1/2}e^{iy t}}{(iy+i\omega)(-y^2-B_+)(-y^2-B_-)}idy.
\end{align*}
Since the integrand in $I_{G'H'}$ is finite, the integral is $O(\varepsilon)$. Therefore, we conclude that $\lim_{\varepsilon\rightarrow 0} I_{G'H'} = 0$. We combine the remaining contributions as
\[
I_{F'I'} = \int_{-R}^{-c\nu_-}\frac{(iy)(-y^2+c^2\nu_-^2)^{1/2}e^{iy t}}{(iy+i\omega)(-y^2-B_+)(-y^2-B_-)}idy.
\]
The integrand is singular at $y=-\omega$ and need further investigation. 

To the integrand, we add and subtract
\[
\frac{iy(-\omega^2+c^2\nu_-^2)^{1/2}ie^{iyt}}{(iy+i\omega)(-y^2-B_+)(-y^2-B_-)},
\]
and regroup the expression as
\begin{align*}
\frac{iy ie^{iyt}[(-y^2+c^2\nu_-^2)^{1/2}-(-\omega^2+c^2\nu_-^2)^{1/2}]}{(iy+i\omega)(-y^2-B_+)(-y^2-B_-)}
+ \frac{iy ie^{iyt}(-\omega^2+c^2\nu_-^2)^{1/2}}{(iy+i\omega)(-y^2-B_+)(-y^2-B_-)}.
\end{align*}
We can easily find the limit of the first expression using l'Hopital's rule and it is finite.

For the second term, we first rewrite it as
\begin{equation*}
 -(\omega^2-c^2\nu_-^2)^{1/2} e^{iyt}  \left[ \frac{y}{(y+\omega)(y^2+B_+)(y^2+B_-)} \right].
\end{equation*}
Next, by performing partial fractions on the rational term, we can further split the expression into 
\begin{align*}
\frac{-(\omega^2-c^2\nu_-^2)^{1/2} e^{iyt}}{(\omega^2+B_+)(\omega^2+B_-)}
&\left[ -\frac{\omega}{(y+\omega)} \right. \\ & + \left. \frac{\omega y^3-\omega^2 y^2 + \omega^3 y + (B_++B_-)\omega y+B_+ B_-}{(y^2+B_+)(y^2+B_-)}\right].
\end{align*}
Along the paths of integration $F'I'$, the second term inside the bracket is finite, so there are no issues in evaluating the integrals. However, the first term in the bracket approaches $+\infty$ as $y\rightarrow -\omega -0$ and $-\infty$ as $y\rightarrow -\omega +0$. The integral of concern is rewritten as
\[
\int_{-R}^{-c\nu_-} \frac{e^{iyt}}{y+\omega} dy = 
\lim_{\varepsilon\rightarrow 0} \int_{-R}^{-\omega-\varepsilon} \frac{e^{iyt}}{y+\omega} dy
+ \lim_{\varepsilon\rightarrow 0} \int_{-\omega+\varepsilon}^{-c\nu_-} \frac{e^{iyt}}{y+\omega} dy.
\]
We split this further as
\begin{equation*}
\mbox{P.V.}   \left[e^{-i\omega t} \int_{-\infty}^\infty \frac{e^{iyt}}{y} dy \right] - \int_{-\infty}^{-R}\frac{e^{iyt}}{y+\omega} dy - \int_{-c\nu_-}^\infty \frac{e^{iyt}}{y+\omega} dy .
\end{equation*}
Evaluating the first term and taking the limit as $R\rightarrow \infty$, we arrive at
\begin{align*}
\lim_{R\rightarrow\infty} I_{F'I'} & =\int_{-\infty}^{-c\nu_-}
\frac{iy ie^{iyt}[(-y^2+c^2\nu_-^2)^{1/2}-(-\omega^2+c^2\nu_-^2)^{1/2}]}{(iy+i\omega)(-y^2-B_+)(-y^2-B_-)} dy \\
& + \alpha \int_{-\infty}^{-c\nu_-}\frac{e^{iyt}(\omega y^3-\omega^2 y^2 + \omega^3 y + (B_++B_-)\omega y+B_+ B_-)}{(y^2+B_+)(y^2+B_-)} dy \\
& + \alpha \left[i\pi e^{-i\omega t} + \int_{-c\nu_-}^\infty \frac{e^{iyt}}{y+\omega} dy \right],
\end{align*}
where $\alpha = -\sqrt{\omega^2-c^2\nu_-^2) }/[(\omega^2+B_+)(\omega^2+B_-)]$.

Putting this all together, we get that the inversion of the Laplace transform yields
\begin{equation}
    \tilde{j}_{b}(t) = \sum_{l=1}^{5}\text{Res}[\hat{\tilde{j}}_b(s) e^{st}, q_l] - \frac{1}{\pi i}\left(  \lim_{R\rightarrow\infty} I_{C'D'}+ \lim_{R\rightarrow\infty} I_{F'I'} \right).
    \label{jb-caseivb}
\end{equation}

\medskip
\noindent Sub-case (c): $c\nu_- = \omega_0 - \omega_1$

In this case, $\hat{\tilde{j}}_{b_2}$ takes the form
\[
\hat{\tilde{j}}_{b_2} = \frac{s(s-ic\nu_-)^{1/2}}{(s+ic\nu_-)^{1/2}(s^2-B_+)(s^2-B_-)}.
\]
We still have branch points at $\pm ic\nu_-$, however, the pole at $-i\omega$ is no longer present. The resulting inverse Laplace transform is
\begin{equation}
\tilde{j}_b(t) = \sum_{l=1}^5 \mathrm{Res}[\hat{\tilde{j}}_{b_1}(s)e^{st}] +
\sum_{l=1}^4 \mathrm{Res}[\hat{\tilde{j}}_{b_2}(s)e^{st}]-\frac{1}{\pi i} \left(
\lim_{R\rightarrow\infty} I_{C'D'} + \lim_{R\rightarrow\infty} I_{F'G'}
\right),
    \label{jb-caseivc}
\end{equation}
where $p_{1,2}=\pm{B_+}$, $p_{3,4}=\pm{B_-}$, and $p_5=-i\omega$. We refer to Figure \ref{bromwich-jb-caseiii}, and
\begin{align*}
I_{C'D'} &= \int_{c\nu_-}^R \frac{y\sqrt{iy-ic\nu_-}e^{iyt}}{\sqrt{iy+ic\nu_-}(y^2+B_+)(y^2+B_-)}dy,\\
I_{F'G'} &= -\int_{-R}^{-c\nu_-} \frac{y\sqrt{iy-ic\nu_-}e^{iyt}}{\sqrt{iy+ic\nu_-}(y^2+B_+)(y^2+B_-)}dy.
\end{align*}

\end{document}